\documentclass[onefignum,onetabnum]{siamart251216}

\usepackage{graphicx}
\ifpdf
  \DeclareGraphicsExtensions{.eps,.pdf,.png,.jpg}
\else
  \DeclareGraphicsExtensions{.eps}
\fi
\usepackage{enumitem}
\usepackage{amsmath, amssymb}
\usepackage{physics}
\usepackage{mathtools}
\usepackage{amsfonts}
\usepackage{titlesec}
\usepackage{subcaption}
\usepackage{algorithm, algpseudocode}

\newcommand{\R}{\mathbb{R}}

\newcommand{\C}{\mathbb{C}}
\newcommand{\innerp}[2]{\left \langle #1, #2 \right \rangle}
\newcommand{\Lind}{L}

\newcommand{\Hil}{\mathcal{H}}

\headers{CPTP Schemes for Lindblad using TT/MPS}{P. DelMastro, D. Appel\"{o}, Y. Cheng}

\title{Completely Positive and Trace Preserving Schemes with Tensor Train Compression for the Lindblad Equation}

\author{
    Peter DelMastro\thanks{Department of Mathematics, Virginia Tech, Blacksburg, VA 24060 (\email{pdelmastro@vt.edu,\\ appelo@vt.edu, yingda@vt.edu}).}
    \and 
    Daniel Appel\"{o}\footnotemark[1]
    \and
    Yingda Cheng\footnotemark[1]
}

\usepackage{amsopn}

\ifpdf
\hypersetup{
  pdftitle={TT CPTP Lindblad},
  pdfauthor={P. DelMastro, D. Appel\"{o}, Y. Cheng}
}
\fi

\begin{document}

\maketitle

%
%
\begin{abstract}
We propose a family of low-rank, completely positive and trace preserving schemes for the Lindblad equation, a common model for open quantum systems.
Low-rank representation is employed at two levels: the density matrix is factorized into the product of tall-skinny matrices, and the columns of these matrices are further represented using the tensor train (TT) format, also know as matrix product states (MPS). 
This two-level low-rank format fits naturally into our existing Kraus is King scheme \cite{Appelo2024-kraus-is-king} for the Lindblad equation, whose underlying operations are arithmetic on the columns of the tall-skinny matrices.
We show how these operations can be performed efficiently in the TT/MPS format, with particular emphasis on density matrix rank-truncation.
We conclude with extensive numerical experiments demonstrating the convergence of this scheme and its efficiency in simulating systems with up to $10^{19}$ degrees of freedom using only modest compute resources.
\end{abstract}

%
%
\begin{keywords}
low-rank methods, tensor trains, Lindblad equation, open quantum systems, completely positive and trace preserving
\end{keywords}

\begin{MSCcodes}
65L99, 
15A69, 
81Q99 
\end{MSCcodes}

%
%
\section{Introduction}

Quantum mechanical phenomena are becoming increasingly relevant in the development of new technologies, with quantum computers being one exciting example.
Numerical simulation of quantum devices is a promising avenue to accelerate this development, but it remains challenging due to the curse of dimensionality.
It is intractable, for instance, to use full quantum mechanical treatments of both the device and its environment, so researchers often derive a reduced model for how a subsystem of interest interacts with its environment. This leads to the study of \textit{open quantum systems} \cite{Breuer2006-Open-Quantum-Systems-Textbook}.

The state of an open quantum system is characterized by its \textit{density matrix} $\rho \in \C^{N \times N}$, where $N$ denotes the size of the Hilbert space of the corresponding closed system. 
This matrix is Hermitian positive (semi-)definite with unit trace, and it encodes the probability that the quantum system is found in any state $x \in \C^N$: $\text{Prob}(x) = x^\dagger \rho x$. 
A density matrix is said to represent a \textit{pure state} when it has rank 1, and it is called a \textit{mixed state} otherwise.

The Lindblad master equation \cite{Lindblad1976-master-eqn, Manzano2020-Lindblad-intro} is one commonly used model for the dynamics of open quantum systems. It describes a system whose interactions with its environment depend only on the system's current state, and it is given by
\begin{equation}
\label{eqn:lindblad}
\dot{\rho} = \mathcal{L}\rho = - i (H\rho - \rho H) + \sum_{j} \left [ \Lind \rho \Lind_j^\dagger - \frac{1}{2} \big( \Lind_j^\dagger \Lind_j \rho + \rho \Lind_j^\dagger \Lind_j \big) \right]. 
\end{equation}
Here, $H = H^\dagger \in \C^{N \times N}$ is the \textit{Hamiltonian}, and $\Lind_j \in \C^{N \times N}$ are called \textit{jump operators}. 

Lindbladian and Schrödinger dynamics are equivalent if all $\Lind_j = 0$ and the system is initially in a pure state.
The inclusion of jump operators results in the deviation from the Schrödinger picture, as these terms in Eqn. (\ref{eqn:lindblad}) drive the system into a mixed state. 
That being said, if the jump operators $\Lind_j$ are small in magnitude compared to the Hamiltonian $H$, the density matrix will often remain \textit{low-rank} for some window of time. 

Past work, e.g. \cite{Appelo2024-kraus-is-king, Chen2025-low-rank-lindblad, cao2025structure, LiWang2023HigherOrder,robin2025unconditionally,le2015adaptive} have taken advantage of this low-rankedness to enable efficient structure-preserving numerical schemes for the Lindblad equation. Most use \textit{Choi's theorem} (c.f. \cite{Manzano2020-Lindblad-intro}), which states that a linear map $\mathcal{G}$ is \textit{complete positivity} (CP) if and only if it can be written in the form
\begin{equation}
    \label{eqn:chois-theorem}
    \mathcal{G}\rho = \sum_{i} G_i \rho G_i^\dagger, \quad G_i \in \mathbb{C}^{N \times N}.
\end{equation}
A CP scheme $\rho^{n+1} = \mathcal{G} \rho^n$ with $\mathcal{G}$ taking this \textit{Kraus form} can be made \textit{trace preserving} (TP) if it is followed by trace normalization: $\rho^{n+1} \gets \rho^{n+1} / \trace{(\rho^{n+1})}$.
Even when the density matrix is low-rank, its representation as $\rho = VV^\dagger$, $V \in \C^{N \times r}$, can still be computationally intractable for larger quantum systems, for which the space $\mathcal{H} = \C^N$ is often the tensor product of smaller spaces: $\mathcal{H} \cong \mathcal{H}_1 \otimes \cdots \otimes \mathcal{H}_d$ with $\mathcal{H}_i = \C^{n_i}$. 
These smaller systems may represent individual magnetic spins or quantum bits (qubits), for which $n_i$ tends to be small.
Still, $N = \textrm{dim}(\mathcal{H}) = \prod_i \textrm{dim}(\mathcal{H}_i) = \prod_i n_i$ grows exponentially in the number of subsystems, making it expensive to represent even a single vector in the full space $\mathcal{H}$.

The tensor product structure of these Hilbert spaces naturally lend themselves to low-rank tensor factorizations of their elements. This is to say that we view elements of $\mathcal{H}$ as order-$d$ tensors $x(i_1, i_2, ..., i_d)$ in $\mathcal{H}_1 \otimes \cdots \otimes \mathcal{H}_d$, and then a compressed tensor format is adopted to represent such tensors.
%
One class of tensor formats, called \textit{tensor networks}, has been particularly prevalent in quantum simulation, with applications including condensed matter physics \cite{White1992-DMRG}, quantum chemistry \cite{Chan2016-DMRG-Review}, and more recently quantum computing \cite{GonzalezGarcia2025-Block-DMRG-X}.
Within tensor networks, tensor trains (TTs) \cite{Oseledets2011-TTs}, also known as matrix product states (MPS) \cite{Schollwock2011-MPS-basics}, have seen the largest success, particularly for eigenvalue calculations \cite{White1992-DMRG, GonzalezGarcia2025-Block-DMRG-X} and time-integration \cite{Schollwock2019-MPS-time-evolution}.

Despite the success of TT/MPS for closed quantum systems, there are, to our knowledge, no schemes for the Lindblad equation that use this tensor format and are also CPTP.
Some researchers use TT/MPS or related formats like \textit{matrix product operators} (MPOs) to represent the density matrix itself \cite{Verstraete2004-MPDOs, Landa2022-LindbladMPO-qiskit, Ceruti2025-TTNOs, Lindoy2025-pyTTN, Houdayer2025-TensorMixedStates, Zhan2026-MPO-Lindblad-ground-state-prep}, but these methods do not maintain positivity.
Others \cite{Werner2014-LPDOs, Yin2024-lpdo-ms-thesis} adopt the positivity-preserving decomposition $\rho = VV^\dagger$ and use the MPO format to represent $V$; this is called the \textit{locally purified density operator} (LPDO) format.
There is also the stochastic approach of quantum trajectories \cite{Dalibard1992-WavefunctionAT-Quantum-Trajectories,Dum1992-MonteCS-QuantumTrajectory}, which can utilize TT/MPS time-integrators, but this technique does not solve the Lindblad equation directly.

Although the LPDO format can efficiently represent open quantum systems, using this format within time-integration schemes comes with the added overhead of \textit{disentangling} to maintain efficient compression \cite{Muller2024-disentangling, Yin2024-lpdo-ms-thesis}.
This optimization procedure, performed at least once per timestep, can be expensive even when the density matrix remains very low rank through time.
In such a case, an alternative is to represent the columns of $V$ individually in TT/MPS format and use their suite of optimization-free compression techniques. 
%
This is the approach we adopt in this paper, extending our recently developed low-rank CPTP scheme \cite{Appelo2024-kraus-is-king} to larger quantum systems.

The rest of this paper is organized as follows. In Section \ref{sec:tt-background}, we review the basics of the TT/MPS format to establish our notation for the paper.
We then review in Section \ref{sec:low-rank-cptp-background} the Kraus-is-King scheme \cite{Appelo2024-kraus-is-king} that we extend to use TT/MPS compression in Section \ref{sec:tt-cptp-scheme}. 
In Section \ref{sec:accelerating-rounding} we discuss methods for accelerating the central truncation routine, and we conclude with numerical examples of the scheme in Section \ref{sec:numerical-experiments}, one on dissipative spin chains and the others pertaining to quantum computing hardware.

%
%
\section{Tensor Train Preliminaries}
\label{sec:tt-background}

\textit{Tensor networks} encompass wide family of tensor formats that have been largely successful in representing quantum systems.
Within this family, \textit{tensor trains} (TTs) \cite{Oseledets2011-TTs, Schollwock2011-MPS-basics}, also called \textit{matrix-product states} (MPS), are particularly well suited to represent systems whose underlying structure is one-dimensional, such as chains of qubits.
This section reviews the ingredients of the TT/MPS format that are necessary for our numerical scheme. 

%
%
%
\subsection{Tensor Trains / Matrix Product States}

The TT/MPS format represents an order-$d$ tensor $x \in \C^{n_1 \times \cdots n_d}$ in terms of a sequence of order-3 tensors $X_k \in \C^{b_{k-1} \times n_k \times b_k}$ called \textit{cores}. Each element of $x$ is expressed as a product of matrices obtained by fixing one element of each core; in particular,
\begin{equation}
x(i_1,...,i_d) = \prod_{k=1}^d X_k(:,i_k,:) = \sum_{\sigma_0=1}^{b_0} \sum_{\sigma_1=1}^{b_1} \cdots \sum_{\sigma_{d}=1}^{b_{d}} \bigg[\prod_{k=1}^d X_k(\sigma_{k-1},i_k,\sigma_k)\bigg].
\end{equation}
The tuple $(b_0,b_1,...,b_d)$ are called the \textit{bond dimensions} or TT-ranks of $x$, where we must have $b_0 = b_d = 1$.
%
Tensor trains are often represented diagrammatically as 
\begin{equation}
    \label{fig:tensor-trains}
    x(i_1,...,i_d) \ = \ \raisebox{-5mm}{\includegraphics[height=1.5cm]{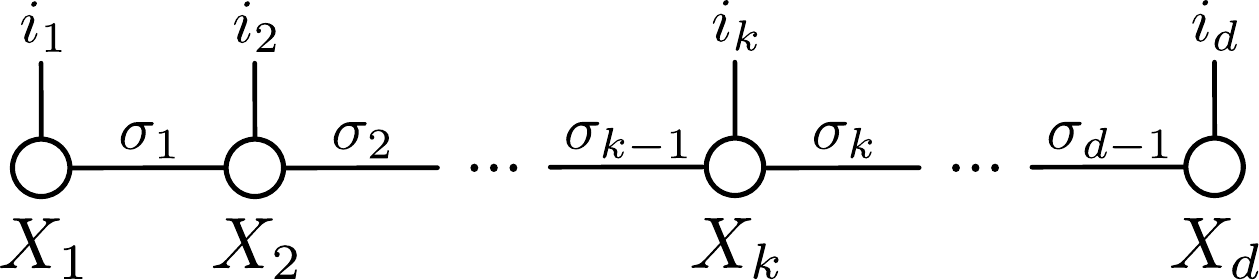}} .
    \end{equation}
The \textit{diagram as a whole} represents the tensor $x(i_1,...,i_d)$ and shows its decomposition as the \textit{contraction} of the cores $X_k$, each represent as a node. Each edge associated with an index, with the \textit{bond indices} $\sigma_1,...,\sigma_{d-1}$ labeling the shared edges.  
The bond dimensions $(b_1,...,b_{d})$ determine the efficiency of this compressed representation: the total degrees of freedom across all the cores is $\sum_{i=1}^d b_{k-1} n_k b_k$.

%
%
\subsection{TT/MPS Addition and Truncation}

Arithmetic operations in the TT format are non-trivial due to the compressed nature of this data structure. Here, the notation $z = x + y$ for two order-$d$ tensors $x$ and $y$ refers to \textit{component wise addition}, namely $z$ is an order-$d$ tensor with elements, i.e. $    z(i_1,...,i_d) = x(i_1,...,i_d) + y(i_1,...,i_d)$. To compute a TT representation of $z = x + y$ where $x$ and $y$ are given in TT format, the standard approach is to first build a sub-optimal TT representation of $z$ and then compress it. Letting $X_k$ and $Y_k$ denote the order-3 cores of $x$ and $y$, respectively, $z$ can be represented in TT format with cores
\begin{align}
    Z_1(i_1,:) &= \begin{bmatrix} X_1(i_1,:) & Y_1(i_1,:)\end{bmatrix},
    \\[0.5em]
    Z_k(:,i_k,:) &= \begin{bmatrix}
        X_k(:,i_k,:) \\
        & Y_k(:,i_k,:)
    \end{bmatrix} 
    \ , \enskip k = 2, ..., d-1,
    \\[0.5em]
    Z_d(:,i_d) &= \begin{bmatrix} X_d(:,i_d) \\ Y_d(:,i_d)\end{bmatrix} .
\end{align}
The resulting representation for $z$ has \textit{increased bond dimensions} $b_k^{(z)} = b_k^{(x)}+b_k^{(y)},$ where $b_k^{(x)}$ and $b_k^{(y)}$ are the bond dimensions of $x$ and $y$, respectively. This representation is nothing more than a relabeling of indices and almost always results in sub-optimal bond dimensions for $z$. 
%
A compression step, often called \textit{truncation}, is therefore employed to reduce the bond dimensions of $z$. There are a few options in this regard, such as the TT-SVD \cite{Oseledets2011-TTs}, randomized algorithms \cite{Daas2021-Randomized-TT-Rounding,Daas2025-AdaptiveRT}, or the TT cross approximation \cite{Oseledets2010-TT-Cross,Savostyanov2011-TT-Cross-Greedy}. 
The TT-SVD is the most accurate truncation algorithm, though also the most expensive. For more details, refer to Section 2.3 of \cite{Daas2021-Randomized-TT-Rounding}, which provides an intuitive presentation of this algorithm by comparing it to a truncation procedure for low-rank matrices.

\subsection{Inner products}

$\langle x, y\rangle = \sum_{i_1,...,i_d} \bar{y}(i_1,...,i_d) x(i_1,...,i_d)$ in the TT format are computed via a series of tensor contractions of the cores of $x$ and $y$. This process is most easily described diagrammatically as
\begin{equation}
    \innerp{x}{y} = \raisebox{-6.5mm}{\includegraphics[height=1.5cm]{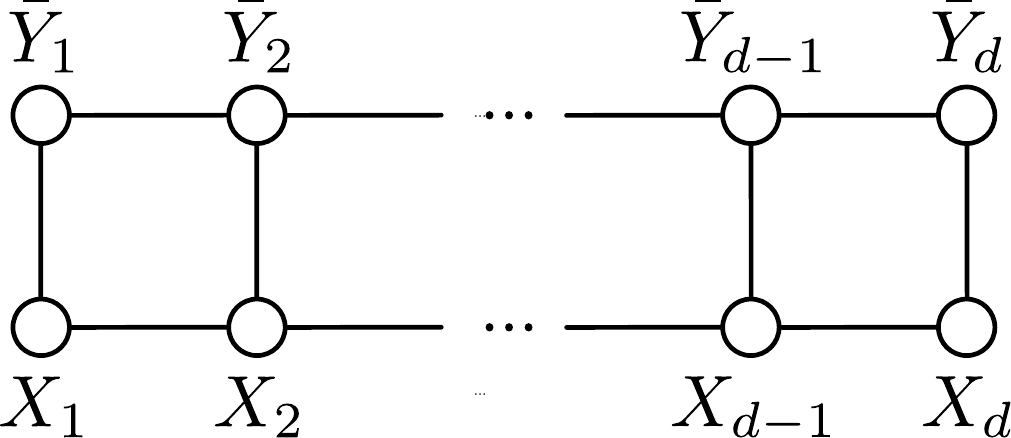}}, 
\end{equation}
where $\bar{Y}_k$ denotes the element-wise complex conjugate of core $k$ of $y$.

\subsection{Matrix Product Operators (MPOs)}

The tensor train format is well-suited to working with operators that can be expressed efficiently in terms of their actions on tensor train cores.
In general, any operator $H$ acting on the Hilbert space $\mathcal{H} \cong \mathcal{H}_1 \otimes \mathcal{H}_2 \otimes \cdots \otimes \mathcal{H}_d$ can be written in \textit{matrix-product operator (MPO) format}, wherein elements of $H$ are expressed as contractions of order-4 tensors $H_k$: 
\begin{equation}
H(i_1,...,i_d ; j_1,...,j_d) 
\ 
= \ \sum_{\gamma_1=0}^{b_0} \sum_{\gamma_2 = 1}^{b_2} \cdots \sum_{\gamma_d = 1}^{b_d} \left[ \prod_{k=1}^{d} H_k(\gamma_{k-1}, i_k, j_k, \gamma_k)\right].
\end{equation}
The number of terms $b_1, b_2, ..., b_d$ in these sums are again called bond dimensions.

As with addition, multiplication of an operator $H$ in MPO format with a tensor $x$ in TT/MPS format involves constructing a sub-optimal TT representation of the product followed by TT truncation.
Letting $b^{(H)}_k$ and $b_k^{(x)}$, $k = 1, 2, ..., d$, denote the bond dimensions of $H$ and $x$, respectively, the tensor $y = Hx$ can be represented as a TT with bond dimensions $b_k^{(y)} = b_k^{(H)} b_k^{(x)}$.
This representation is achieved by setting the cores $Y_k$ of $y$ to be the contraction of cores $H_k$ and $X_k$ of $H$ and $x$. See \cite{Schollwock2011-MPS-basics} for more details. The sub-optimal representation of $y = Hx$ is followed by truncation to lower bond dimensions.
This compression needed especially for MPOs with large bond dimension due to the multiplicative increase in the MPS bond dimension.
Randomized rounding \cite{Daas2021-Randomized-TT-Rounding, Daas2025-AdaptiveRT} can be particularly effective in this regard because these methods operate without explicitly forming the sub-optimal MPS representation of $y$.

%
%
\section{Kraus is King: A low-rank CPTP scheme for the Lindblad equation}
\label{sec:low-rank-cptp-background}

We recently developed a CPTP scheme for the Lindblad equation by introducing an integrating factor to handle the terms that aren't in Kraus form \cite{Appelo2024-kraus-is-king}.
The current paper directly extends our previous work by introducing TT/MPS compression within the low-rank factorization of the density matrix.
We therefore begin by reviewing the primary steps of this ``Kraus is King'' method for the Lindblad equation.
Algorithm (\ref{alg:kraus-is-king}) summarizes a single timestep of the scheme, and the rest of this section describes the algorithm in more detail.

%
%
\begin{algorithm}[t!]
\caption{Low-rank CPTP time-integrator for the Lindblad equation}
\label{alg:kraus-is-king}
\begin{algorithmic}[1]
    \Require Factor matrix $V \in \C^{N \times r}$ of $\rho = VV^\dagger = \rho(t)$;  Butcher Tableau $\mathcal{S} = (\mathbf{A}, \mathbf{b}, \mathbf{c})$ of an $s$-stage explicit scheme; \ Step size $h $; \  Truncation tolerance $\tau$.
    \vspace{0.35\baselineskip}
    \Ensure Updated factor $V' \in \C^{N \times r'}$ of $\rho(t+h)$.
    \vspace{0.35\baselineskip}
    \State $\tau \ \gets \  \tau \ / \ (s+1)$
    \Comment{Tolerance per compression/truncation within the scheme} 
    \vspace{0.35\baselineskip}
    \Statex \textbf{Intermediate Stages}
    \State $V^0 := V$
    \For {$i = 1, 2, ..., s$}
        \State $U^{i} = \texttt{Schrodinger-Solve}(V, c_i h)$ \Comment{Solve $\dot{\Psi} = -iH_\textrm{eff}\Psi$ with initial state $V$}
        \State \Comment{and time step $c_i h$ \hspace{3.2cm} \ } 
        \State $W^{i-1} = \left[L_1 V^{i-1}, L_2 V^{i-1}, ..., L_{P} V^{i-1}\right]$ 
        \Comment{Applying jump operators $L_p$}
        \For {$j = 1, 2, ..., i-1$}
            \State $Y^{i,j} = \texttt{Schrodinger-Solve}(\sqrt{a_{ij} h} \ W^j, \ (c_i-c_j) h)$
        \EndFor
        \State $V^{i} = \texttt{Compress}\left(\left[U^i, \ Y^{i,1}, ...,  \ Y^{i,i-1} \right], \ \tau \right)$
	\EndFor 
    \vspace{0.35\baselineskip}
    \Statex \textbf{Final Stage}
    \State $U = \texttt{Schrodinger-Solve}(V, h)$
    \State $W^{s} = \left[L_1 V^{s}, L_2 V^{s}, ..., L_{P} V^{s}\right]$ 
    \For {$i = 1, 2, ..., s$}
        \State $Y^{i} = \texttt{Schrodinger-Solve}(\sqrt{b_i h} \ W^i, (1-c_i) \ h)$
    \EndFor
    \State $V' = \texttt{Compress}\left(\left[U,  \ Y^{1}, ...,  \ Y^{s}\right], \ \tau \right)$ 
    \vspace{0.35\baselineskip}
    \State $V' \gets V' / \norm{V'}_F$ \Comment{Trace normalization}
\end{algorithmic}
\end{algorithm}

Given a Butcher Tableau ($\mathbf{A}, \mathbf{b}, \mathbf{c}$) defining an $s$-stage method, the Kraus is King method computes the next state $\rho^{n+1} = \rho(t_n+h)$ from $\rho^n = \rho(t_n)$ in stages as
\begin{align}
   \rho^{n,i} &= U(c_i h) \rho^n U(c_i h)^\dagger + h \sum_{j < i} \ a_{i,j} \ U((c_i-c_j) h) \ \mathcal{L}_L \rho^{n,j} \ U((c_i-c_j) h)^\dagger,
   \\[0.5em]
   \rho^{n+1} &= U(h) \rho^{n} U(h)^\dagger + h \sum_{i=1}^s \ b_i \ U((1-c_i)h) \ \mathcal{L}_L \rho^{(n,i)} \ U((1-c_i)h)^\dagger.
\end{align}
Here, $U(ch) \approx \exp(-c h H_\textrm{eff})$ denotes an approximate flow operator with respect to the effective Hamiltonian
\begin{equation}
    \label{eqn:H-eff}
     H_\textrm{eff} = H + \frac{1}{2i} \sum_{j} \Lind_j^\dagger L_j.
\end{equation}
The density matrices $\rho^n$ are further factorization as $\rho^n = (V^n)(V^n)^\dagger$ where $V^n \in \C^{N \times r_n}$ has only as many columns as the rank of $\rho^n$.
The intermediate stages $\rho^{n,i}$ are similarly stored only in this low-rank format, with a \textit{compression} step introduced to maintain optimal representation. Setting $\hat{h}_i = c_i h$, $h_{i,j} = h_i - h_j$, and $\alpha_{i,j} = \sqrt{a_{ij} h}$ to simplify notation, the low-rank stage calculations are performed as
\begin{align}
    \tilde{V}^{n,i} &= \left[ U(\hat{h}_i) V^n, \ \alpha_{i,1} \ U(h_{i,1}) \ \mathcal{L}_L V^{n,1}, \ \cdots, \ \alpha_{i,i-1} \ U(h_{i,i-1}\big) \mathcal{L}_L V^{n,i-1} \right] ,
    \\
    V_n^{i} &= \texttt{Compress}(\tilde{V}_n^{(i)}, \tau) .
\end{align}
Above, $\texttt{Compress}(\cdot, \tau)$ is a CP truncation function that, given input $X \in \C^{N \times r}$, returns a new factor matrix $\tilde{X} \in \C^{N \times \tilde{r}}$ with $\tilde{r} \le r$ such that $\lVert XX^\dagger - \tilde{X}\tilde{X}^\dagger \rVert_F \le \tau$. To ensure convergence of the method, the truncation tolerance for the intermediate states should be taken as $\tau = O(h^{p+1/2})$ where the local truncation error of the RK scheme is $O(h^{p+1})$.  The tolerance when computing $V^{n+1}$ should be $O(h^{p+1})$. The latter choice follows from standard local vs. global truncation error analysis, and for the former, the half power is due to working with Cholesky factors. 

The truncation function $\texttt{Compress}(\cdot, \tau)$ in our previous paper was computed via an SVD and proceeds as follows.
Given a matrix $X \in \C^{N \times r}$, first orthogonalize as $X = QR$, and then compute $R = U\Sigma V^\dagger$. These operations yield an SVD factorization $X = (QU)\Sigma V^\dagger$ of $X$, from which one sets $\tilde{X} = Q U(:,:\tilde{r})\Sigma(:\tilde{r},:\tilde{r})$, where the truncation rank $\tilde{r}$ is selected as the smallest integer such that
\begin{equation}
    \norm{XX^\dagger - \tilde{X}\tilde{X}^\dagger}_F = \sum_{i=\tilde{r}+1}^{r} \sigma_i^2 \le \tau . 
\end{equation}
This map can be shown to be CP, c.f. Section 2 of \cite{Appelo2024-kraus-is-king}.

%
%
\section{CPTP Scheme with Tensor Train Compression}
\label{sec:tt-cptp-scheme}

We now present our extension of the Kraus-is-King scheme \cite{Appelo2024-kraus-is-king} to use the TT/MPS format to further compress the factor $V$ of $\rho = VV^\dagger$. 
We begin in Section \ref{sec:tt-format-for-density-matrices} by discussing at a high-level our use of this tensor format and how this choice necessitates changes within the Kraus-is-King scheme.
We then describe the main subroutine -- density matrix rank compression -- with our selected tensor format.

\subsection{TT/MPS Compressed Low-Rank Density Matrix Format}
\label{sec:tt-format-for-density-matrices}

\begin{figure}
    \centering
    \includegraphics[width=0.7\linewidth]{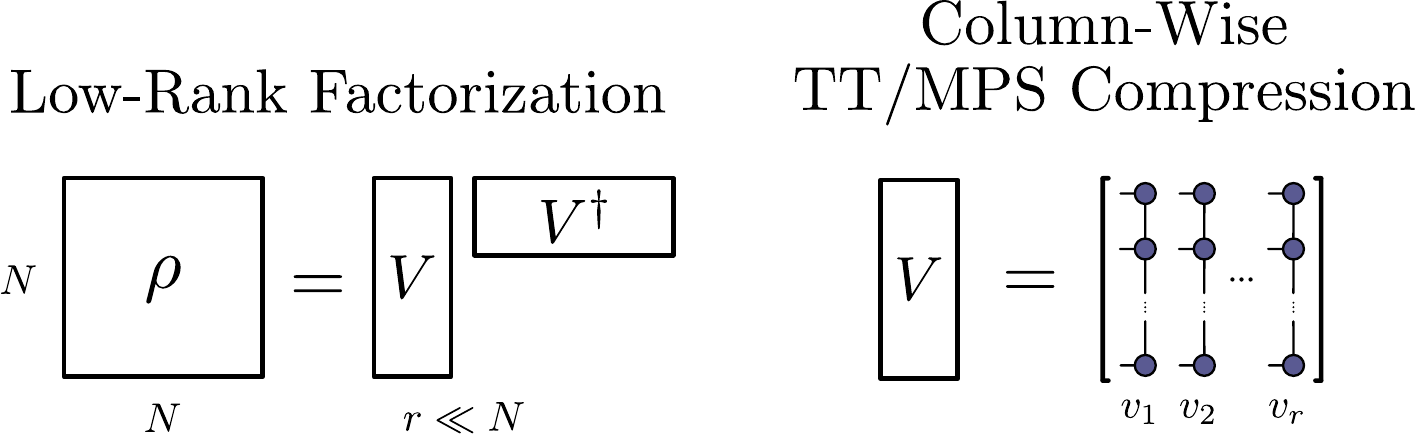}
    \caption{Two-level low-rank decomposition of the density matrix $\rho \in \C^{N\times N}$. At the top level, $\rho$ is factorized as $\rho = VV^\dagger$ where $V \in \C^{N \times r}$, $r \ll N$, is a tall-skinny matrix. At the second level, the columns $v_j \in \C^N$ of $V$ are viewed as order-$d$ tensors of shape $n_1 \times \cdots \times n_d$ with $N = n_1 n_2 \cdots n_d$, and each $v_j$ is represented in TT/MPS format with its own cores $G_{j,k}$. 
    }
    \label{fig:two-level-low-rank-format}
\end{figure}

Let the overall Hilbert space $\mathcal{H} = \C^N$ have tensor product structure $\mathcal{H} = \Hil_1 \otimes \cdots \otimes \Hil_d$, where each $\Hil_k \cong \C^{n_k}$ and $N = \prod_{k=1}^d n_k$. Then each $\psi \in \Hil$ can be viewed as an order-$d$ tensor with elements $\psi(i_1,...,i_d)$ where index $i_k \in \{1,2,...,n_k\}$ can take $n_k$ values.

For systems with this tensor-product structure, we propose to represent the density matrix $\rho \in \C^{N \times N}$ with two levels of low-rankedness, as shown in Figure \ref{fig:two-level-low-rank-format}.
At the top level, we follow the Kraus-is-King scheme and represent the density matrix as $\rho = VV^\dagger$ where $V \in \C^{N \times r}$ is a tall skinny matrix.
As a second level of compression, the columns $v_j \in \C^N$ of $V$ themselves, viewed as order-$d$ tensors, are then represented in the tensor train format, e.g.
\begin{equation}
    v_j(i_1,...,i_d) = G_{j,1}(i_1,:) G_{j,2}(:,i_2,:) \cdots G_{j,d}(:,i_d) , 
\end{equation}
with $G_{j,k} \in \C^{b_{j,k} \times n_k \times b_{j,k+1}}$, for some sequences of cores $G_{j,k}$. 
Each column $v_j$ of $V$ has its own set of cores $\{G_{j,k}\}_{k=1}^d$ and bond dimensions $b_{j,k}$.
Our format for $V$ is distinct from the block TT format, where all vectors share all but one core, often used to solve eigenvalue problems with TTs \cite{Dolgov2013-BlockTT, Kressner2014-BlockTT}.

The tensor format proposed above fits naturally into the Kraus-is-King scheme, which relies only on two types of arithmetic operations on the columns of $V$: taking linear combinations and applying linear operators.
Both of these operations can be performed in the TT/MPS format provided the Hamiltonian and jump operators are amenable to representation as MPOs.

Nearly identical to the original scheme shown in Algorithm (\ref{alg:kraus-is-king}), our new scheme, Algorithm (\ref{alg:tt-kraus-is-king}), indicates the regions needing tensor train arithmetic by using the prefix \texttt{TT-} for the subroutines.
In particular, TT/MPS operations come into play at three points during each stage: (1) solving the Schrödinger equation, (2) applying the jump operators, (3) truncating the density matrix. 
The next few sub-sections discuss these subroutines of the scheme in further detail. 
We elect to start from rank truncation as this routine yields the dominating cost of the scheme.

%
%
\begin{algorithm}[t!]
\caption{Kraus is King integrator with tensor train (TT) compression}
\label{alg:tt-kraus-is-king}
\begin{algorithmic}[1]
\Require Factor matrix $V \in \C^{N \times r}$ with columns $v_i$ represented in TT format;  Butcher Tableau $\mathcal{S} = (\mathbf{A}, \mathbf{b}, \mathbf{c})$; \ Step size $h $; \  Truncation tolerance $\epsilon$.
\vspace{0.35\baselineskip}
\Ensure Factor $V' \in \C^{N \times r'}$ of $\rho(t+h)$
\vspace{0.35\baselineskip}
\State $\tau \ \gets \  \tau \ / \ \left(3(s+1)\right)$
\Comment{Error tolerance per truncation within the scheme}
\vspace{0.35\baselineskip}
\Statex \textbf{Intermediate Stages}
\State $V^0 := V$
\For {$i = 1, 2, ..., s$}
    \State $U^{i} = \texttt{TT-Schrodinger-Solve}(V, \ c_i h, \  \tau)$ 
    \State $W^i = \texttt{TT-Compress-L}\left(\left[L_1 V^{i-1}, ..., L_{P} V^{i-1}\right], \ \tau\right)$ 
    \For {$j = 1, 2, ..., i-1$}
        \State $Y^{ij} = \texttt{TT-Schrodinger-Solve}( \sqrt{a_{ij} h} \ W^j, \ (c_i-c_j) h, \ \tau/i)$
    \EndFor
    \State $V^{i} = \texttt{TT-Compress}\left(\left[U^i, Y^{i,1}, ..., \ Y^{i,i-1} \right], \ \tau\right)$
\EndFor
\vspace{0.35\baselineskip}
\Statex \textbf{Final Stage}
\State $U = \texttt{TT-Schrodinger-Solve}(V, h, \ \tau)$
\State $W^s = \texttt{TT-Compress-L}(\left[L_1 V^{s}, L_2 V^{s}, ..., L_{P} V^{s}\right], \ \tau)$ 
\For {$i = 1, 2, ..., s$}
\State $Y^{i} = \texttt{TT-Schrodinger-Solve}(\sqrt{b_{i}h} \ W^i, \ (1-c_i) h, \ \tau/s)$
\EndFor
\State $V' = \texttt{TT-Compress}\left(\left[U, \ Y^{1}, ...,   \ Y^{s}\right], \ \tau \right)$
\vspace{0.35\baselineskip}
\State $V' \gets V' / \norm{V'}_F$ \Comment{Trace normalization}
\end{algorithmic}
\end{algorithm}

%
%
%
\subsection{\texttt{TT-Compress}: Density matrix rank truncation}
\label{sec:rank-truncation} 

In this section, we introduce the core components used in Algorithm \ref{alg:density-matrix-truncation}.

Let $\rho = XX^\dagger$ be a density matrix where the factor matrix $X = [x_1,...,x_R] \in \C^{N \times R}$ has columns $x_i$ represented in tensor train format. We aim to find a smaller factor matrix $\tilde{X} \in \C^{N \times r}$ such that $\tilde{\rho} = \tilde{X}\tilde{X}^\dagger$ satisfies $\norm{\rho - \tilde{\rho}}_F \le \tau$.

As discussed in Section \ref{sec:low-rank-cptp-background}, the standard approach in dense arithmetic would be to first compute a QR decomposition $X = QR$. With the columns $x_i$ represented in the TT format, performing this decomposition would be particularly costly because all arithmetic operations must also be performed in TT format. The resulting matrix $Q$ would have also have columns $q_i$ in TT format, and their bond dimensions would likely be much higher than those of each $x_i$ in order to make $Q$ close to orthogonal. 

An alternative approach is to compute the eigenvalue decomposition of the inner product matrix $X^\dagger X$, similar to the idea of Cholesky QR, as this gives us the eigenvalues and right singular vectors of $X$: $X^\dagger X = (U\Sigma V^\dagger)^\dagger(U\Sigma V^\dagger) = V \Sigma^2 V^\dagger$. In doing so, we only need to compute $O(R^2)$ pairwise inner products $\langle x_i, x_j \rangle$ in TT format, which are generally cheaper than the $O(R^2)$  TT sums needed for a Gram-Schmidt QR, for instance (here we advise the reader of the clash of notation with $R$ denoting the number of columns in $X$ and the $R$ matrix from the QR factorization). 
If the condition number of $X^\dagger X$ is of concern, yet another approach is by to perform a column-pivoted Cholesky factorization of $X^\dagger X = (LP)^\dagger(LP)$ and then compute $LP = U \Sigma V^\dagger$.

With $\Sigma$ and $V$ computed, the truncation process proceeds to select the smallest index $r$ such that the truncation error satisfies $\sum_{i=r+1}^R \sigma_i^2 \le \alpha \ \tau$. Here, $\alpha \in [0,1]$ is a parameter that controls the proportion of our overall truncation error budget $\tau$ that we allot to error due to this SVD truncation. TT arithmetic will be necessary during the next step of the algorithm, so taking $\alpha < 1$ provides room for this arithmetic to be inexact. Note also that the threshold is $\tau$ here rather than $\tau^2$, as $\sigma_i$ are the singular values of $X$, whereas we want to truncate $\rho = XX^\dagger$ to tolerance $\tau$. 

The SVD-optimal matrix $\tilde{X}^\textrm{exact} = U_r \Sigma_r $ with $V_r = V(:,:r)$, for instance, can be recovered as $\tilde{X}^\textrm{exact} = X V_r$. This is to say that the columns $\tilde{x}_i^\textrm{exact}$ of $\tilde{X}^\textrm{exact}$ are linear combinations of the columns of $X$ as
\begin{equation}
    \tilde{x}_i^\textrm{exact} = \sum_{j=1}^R V(j,i) x_j . 
\end{equation}
These linear combinations are computed \textit{inexactly} in the TT format, either via iterative TT-SVD or randomized TT rounding, which results in the final factor matrix $\tilde{X} = \tilde{X}^\textrm{exact} + E$. The truncation error $E = [e_1, e_2, ..., e_r]$ on $\tilde{X}^\textrm{exact}$ results in  error at the level of the density matrix as
\begin{equation}
    \label{eqn:tt}
    \norm{\tilde{X}\tilde{X}^\dagger - \tilde{X}^\textrm{exact}(\tilde{X}^\textrm{exact})^\dagger}_F
    \le 
    \sum_{i=1}^r \tau_i , \quad  \tau_i : = \norm{e_i}_2^2 + 2 \norm{e_i}_2 \sigma_i . 
\end{equation}
The overall error truncation error on $\rho$ is then
\begin{equation}
    \label{eqn:trunc-error-bound}
    \norm{\rho - \tilde{\rho}}_F \le \underbrace{\sum_{i=r+1}^R \sigma_i^2}_{\textrm{SVD truncation of} \ \rho} + \underbrace{\sum_{i=1}^r \left(\norm{e_i}_2^2 + 2 \norm{e_i}_2  \sigma_i \right)}_{\textrm{TT arithmetic truncation error}} . 
\end{equation}
Determining the optimal balance between errors arising from the SVD truncation and inexact TT arithmetic is generally not feasible, as we do not know how much error will be incurred during the linear combinations without first performing them. Heuristic approaches must therefore be employed to ensure the sum of these terms does not exceed the tolerance $\tau$. We will discuss this trade-off, as well as other accelerations techniques in Section \ref{sec:accelerating-rounding}.

\subsection{Truncation as a CP Map}
\label{sec:truncation-is-CP}

The truncation routine takes an input $XX^\dagger$ with $X \in \C^{N \times R}$ and returns $\tilde{X} \tilde{X}^\dagger$ with $\tilde{X} \in \C^{N \times r}$, $r \le R$. Although our procedure is based on the SVD-based scheme that we previously showed to be a completely positive (CP) map \cite{Appelo2024-kraus-is-king}, the new method introduces TT/MPS truncation when computing the columns of $\tilde{X}$. Here we argue that the truncation remains CP despite these truncations.

Consider the operator $G = \tilde{X} D X^+$ where $X^+$ denotes a pseudo-inverse of $X$, and $D \in \C^{r \times R}$ is a diagonal matrix 
\begin{equation}
    D = \begin{bmatrix}
    1 &   &       &     & 0      & \cdots & 0 \\
      & 1 &       &     & 0      & \cdots & 0\\
      &   & \ddots &    & \vdots & \dots  & \vdots \\
      &   &        & 1  & 0      & \cdots & 0
    \end{bmatrix} 
    = 
    \begin{bmatrix}
        I_{r\times r} & 0_{r\times (R-r)} 
    \end{bmatrix}  .
\end{equation}
Above, $I_{r\times r}$ is an $r \times r$ identity matrix and $0_{r \times (R-r)}$ denotes an $r \times (R-r)$ matrix with all zeros. The action of $G$ on $X$ is then
\begin{equation}
    G X \ = \ (\tilde{X} D X^\dagger) X = \tilde{X} \begin{bmatrix}
        I_{r \times r} & 0_{r \times (R-r)} 
    \end{bmatrix} 
    I_{R\times R}
    \ = \ \begin{bmatrix} \tilde{X} & 0_{N \times (R-r)} \end{bmatrix} , 
\end{equation}
and hence
\begin{align*}
    G(XX^\dagger) G^\dagger = (GX)(GX)^\dagger &= 
    \begin{bmatrix} \tilde{X} & 0_{N \times (R-r)} \end{bmatrix}
    \begin{bmatrix} \tilde{X}^\dagger \\ 0_{(R-r) \times N} \end{bmatrix} 
    \\[0.5em]
    &=\tilde{X}\tilde{X}^\dagger + 0_{N \times (R-r)} 0_{(R-r) \times N} .
\end{align*}
The function $\rho \to G \rho G^\dagger$ with $G$ as above therefore maps $XX^\dagger$ to $\tilde{X}\tilde{X}^\dagger$. It has also been expressed in Kraus form, so the map is CP.

%
%
%
\subsection{\texttt{TT-Compress-L}: Applying jump operators}
\label{sec:jump-operator-compression}

For many systems of interest, the jump operators $L_j$ encoded by operator $\mathcal{L}_L$ individually act on a single subsystem, e.g. a single spin or single qubit.
These operators have Kronecker product structure $L_j = I^{\otimes (j-1)} \otimes \hat{L}_j \otimes I^{\otimes (d-j)}$ where $\hat{L}_j \in \C^{n_{k_j} \times n_{k_j}}$ acts linearly on the $k_j$-th subsystem $\Hil_{k_j} \cong \C^{n_{k_j}}$.
Applying such an operator to a $v \in \Hil$ represented in TT format amounts to contracting this small matrix $\hat{L}_j$ with core at site $k_j$ of $v$. The cores $W_k$ of the resulting tensor $w = L_j v$ will share most of the cores $V_k$ of $v$, aside from core $k_j$:
\begin{equation}
    W_{k_j}(:,i_k,:) = \sum_{\ell} \hat{L}_j(i_k,\ell) V_{k_j}(:,\ell,:)  \quad \textrm{and} \quad W_{k} = V_k  \ \textrm{for} \ k \neq k_j .
\end{equation}
Computing $\mathcal{L}_L V = [\mathcal{L}_L v_1, ..., \mathcal{L}_L v_r]$ for such jump operators is not all too expensive therefore as only $dr$ distinct small contractions must be performed. At the same time, the full matrix $\mathcal{L}_L V$ will have $d$ times as many columns as $V$, 
\begin{equation}
    \label{eqn:LV-full}
    \mathcal{L}_L V = [L_1 v_1, ..., L_d v_1, L_1 v_2, ..., L_1 v_r, ..., L_d v_r] \in \C^{N \times dr} .
\end{equation}
We find that the matrices $\mathcal{L}_L V$ often have many small singular values that can be discarded using the rounding algorithm just outlined in Section \ref{sec:rank-truncation}.
When performing this truncation, one can take advantage of the \textit{shared-core structure} of the columns of $\mathcal{L}_L v =  [L_j v]_{j=1}^d \in \C^{N \times d}$ to improve the efficiency of inner product calculations and linear combinations within the truncation routine. 

\subsubsection{Pairwise inner products} 

$\langle L_j v, L_k v \rangle$ for $k,j=1,...,d$ can be computed using only $O(d^2)$ contractions, whereas $O(d^3)$ contractions would be needed for this many generic tensors. Letting $U_1$, $U_2$, ..., $U_{k-1}$, $C_k$, $V_{k+1}$, ..., $V_{d}$ denote the cores of $v$ in its $k$-th canonical form \cite{Schollwock2011-MPS-basics}, $\langle L_j v, L_k v \rangle$ can be computed by contracting only $|j-k|+1$ cores as
\begin{equation}
    \langle L_j v, L_k v \rangle = \enskip 
	\raisebox{-9.mm}{\includegraphics[height=2cm]{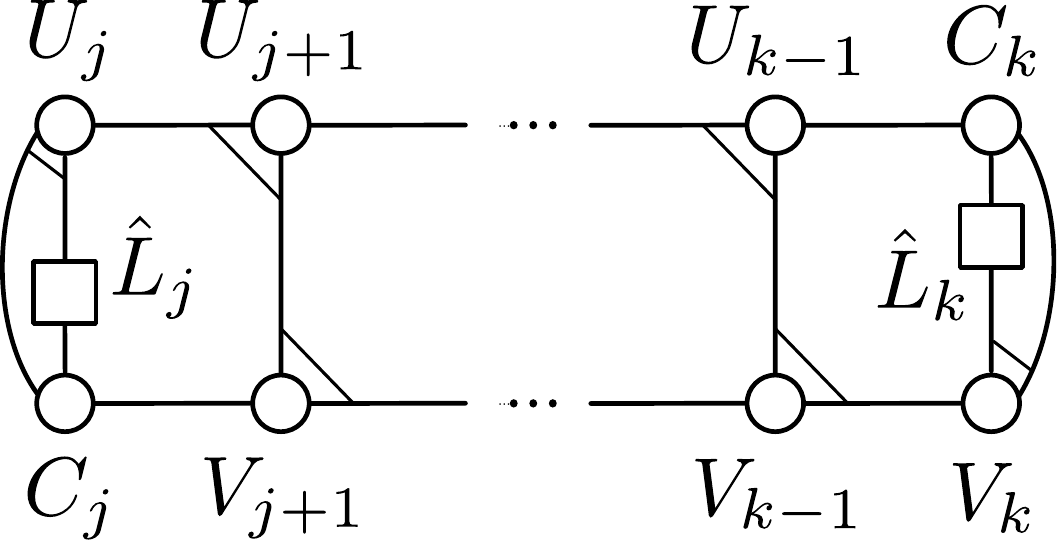}} \ .
\end{equation}
Moreover, the partial contractions involved in computing $\langle L_j v, L_k v \rangle$ can be re-used when computing $\langle L_j v, L_\ell v \rangle$ for $\ell > k$. Let $P_{jk}$ denote the order-2 tensor formed by contracting cores $1$ through $k$ of $L_j v$ and $L_k v$, e.g.
\begin{equation}
    P_{jk} = \raisebox{-9.mm}{\includegraphics[height=2cm]{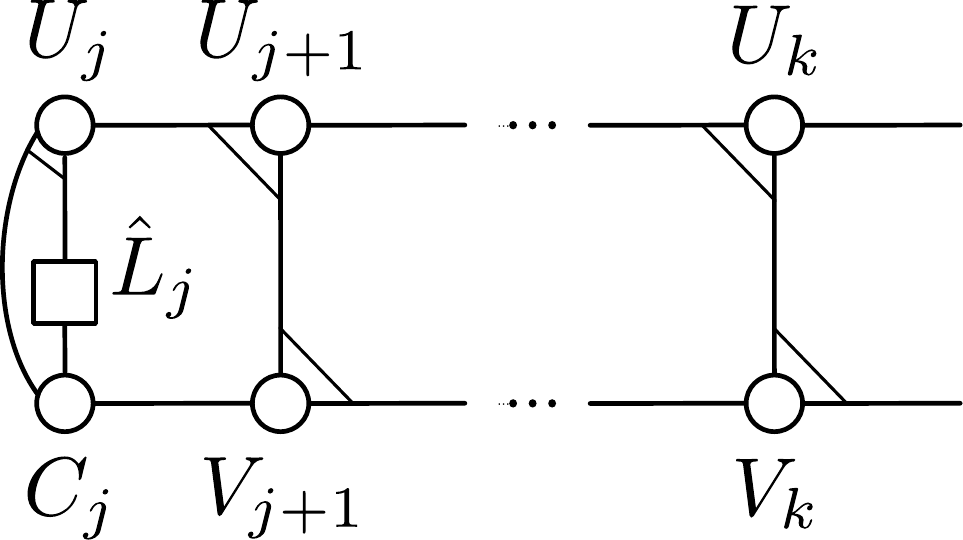}}  . 
\end{equation}
The inner product $\langle L_j v, L_{k+1} v \rangle$ is formed by contracting $P_{jk}$ with $C_k$, $\hat{L}_k$ and $V_{k}$. 
Similarly, the next partial overlap $P_{j,k+1}$, which is used to compute $\langle L_j v, L_{k+2} v \rangle$, is formed by contracting $P_{jk}$ with $U_{k+1}$ and $V_{k+1}$. 
For a fixed $j$, this recurrence allows us to compute all $\langle L_j v, L_{k} v \rangle$ for $k > j$ using only $O(d-j)$ contractions, so only $O(d^2)$ contractions are needed to compute all $\langle L_j v, L_{k} v \rangle$ overall. For generic tensors $x_1, ..., x_d$, each $\langle x_j, x_k \rangle$ needs $d-1$ contractions, so $O(d^3)$ contractions would be needed to compute all of their pairwise inner products.

\subsubsection{Linear combinations} 

$y = \sum_{j=1}^d c_j L_j v$ can be represented as a TT tensor with bond dimensions at most twice those of $v$. Again letting $U_1$, $U_2$, ..., $U_{k-1}$, $C_k$, $V_{k+1}$, ..., $V_{d}$ denote the cores of $v$ in its $k$-th canonical form, the linear combination $y$ can be written in TT format with cores
\begin{align}
    Y_1(i_1,:) &= \begin{bmatrix} W_1(i_1,:) & U_1(i_1,:)\end{bmatrix} ,
    \\[0.5em]
    Y_k(:,i_k,:) &= \begin{bmatrix}
        V_k(:,i_k,:) \\
       W_k(:,i_k,:) & U_k(:,i_k,:)
    \end{bmatrix} 
    , \enskip k = 2, ..., d-1,
    \\[0.5em]
    Y_d(:,i_d) &= \begin{bmatrix} V_d(:,i_d) \\ W_d(:,i_d)\end{bmatrix} \ , 
    \\[0.5em]
    W_k(:,i_k,:) &= \sum_{j \ \textrm{s.t.} \ k_j = k} c_{j} \left[ \sum_{\ell} \hat{L}_j(i_k,:) C_{k_j}(:,i_k,:)\right] .
\end{align}

\subsubsection{Within the CPTP scheme} we make use of these more efficient arithmetic operations to compress each submatrix $\mathcal{L}_L v_i$ of the large matrix $\mathcal{L}_L V$. Symbolically, we take $\tilde{W}_i = \texttt{TT-Compress}\left[ L_1 v_i, ..., L_d v_i\right]$ where this truncation is performed using the specialized inner product and linear combinations routines described above. We then aggregate these individually compressed submatrices into $ W = \texttt{TT-Compress}\big[\tilde{W}_1, ..., \tilde{W}_r\big]$ as a compressed representation of $\mathcal{L}_L V$.

%
%
%
\subsection{Solving the Schrödinger equation in TT/MPS format}
\label{sec:TT-time-evolution}

There are numerous methods for solving the Schrödinger equation in the TT/MPS format, c.f. the review \cite{Schollwock2019-MPS-time-evolution}.
Within our scheme, we use methods that build explicit MPO approximations of the flow operator $U_{h} = \exp\{-ih H_\textrm{eff}\}$, with the effective Hamiltonian $H_\textrm{eff}$ as in Eqn. \ref{eqn:H-eff}. 
Once $U_{h}^\textrm{MPO}$ constructed, flow with respect to the Hamiltonian simply corresponds to MPO-MPS multiplication followed by a compression, either with TT-SVD or randomized rounding.

The operator splitting technique, commonly called \textit{time-evolving block decimation} (TEBD) in the tensor network community \cite{Vidal2004-TEBD, Verstraete2004-TEBD}, is most suitable for systems with nearest neighbor interactions between dimensions. See Supplemental Material \ref{SM-sec:TEBD}   for more details. For systems with longer-range interactions but low bond dimension, the truncated Taylor series approximation offers a simple approach to approximating the flow operator. This is constructed iteratively, with MPO compression performed at each stage to maintain low bond dimension. Starting from $U^{(0)} = \delta U^{(0)} = I$, 
\begin{align}
    U^{(k)} 
    &= 
    \texttt{MPO-Compress}\left(U^{(k-1)} + \delta U^{(k)}, \ \tau \right),
    \\
    \delta U^{(k)}
    &= \texttt{MPO-Compress}\left(-\frac{i h H_\textrm{eff}}{k} \delta U^{(k-1)}, \ {\tau}\right).
\end{align}
This approach to approximating the flow operator is only efficient so long as $U^{(k)}$ and $\delta U^{(k)}$ can be computed with small truncation tolerance $\tau$ without inducing large growth in bond dimensions.
Whether this result is observed depends on the Hamiltonian $H$, but we found this approach useful in of our examples (c.f. Section \ref{sec:mock-quantum-circuit-sim}).

%
%
%
\subsection{Trace Normalization}

The final step of Algorithm \ref{alg:tt-kraus-is-king} is to re-normalize the trace of the density matrix.
With $\rho = VV^\dagger$, $\trace(\rho) = \trace(VV^\dagger) = \norm{V}_F^2$, so trace normalization can be achieved by setting $V \gets V / \norm{V}_F$. 
As we have represented the columns $v_i$ of $V = [v_1, ..., v_r]$ in TT/MPS format, we compute this norm as $\norm{V}_F^2 = \sum_{i=1}^r \norm{v_i}_F^2$.
With careful bookkeeping, the norms $\norm{v_i}_F$ can be determined within the previous \texttt{TT-Compress} step with minimal added cost because the compression step will generates $v_i$ that are \textit{canonicalized}. In such a case, $\norm{v_i}_F$ is simply the Frobenius norm of the TT/MPS core of $v_i$ at its center of orthogonality, typically the leftmost or rightmost site depending on the implementation.

%
%
\section{Details for the Density Matrix Compression}
\label{sec:accelerating-rounding}

Linear combinations in the tensor train format are the primary bottleneck of the truncation algorithm described in Section \ref{sec:rank-truncation}. Here we report a few techniques we found to be effective in reducing the number of TT sums in these calculations while only introducing error on the order of $O(h^{p+1})$. Following the notation of Section \ref{sec:rank-truncation}, we aim to truncate the matrix $\rho = XX^\dagger$ to $\tilde{\rho} = \tilde{X} \tilde{X}^\dagger$, where $X$ has $R$ columns and $\tilde{X}$ has $r < R$ columns, such that $\norm{\rho - \tilde{\rho}}_F \le \tau $.

\subsection{Norm Screening} 

For systems with dissipation, we find many columns $x_i$ of $X$ have small norm, well below the truncation threshold $\tau$. Some  of these columns can therefore be discarded from the density matrix without us needing to compute inner products $x_i^* x_j$ nor include them in linear combinations when forming each $\tilde{x}_i$. Letting $X = [x_1, ..., x_R]$ have columns ordered such that $\norm{x_i}_F \le \norm{x_{i+1}}_{F}$, we select the largest $R_0 \le R$ such that 
\begin{equation}
    \label{eqn:norm-screening}
    \norm{XX^\dagger - X_{R_0} X_{R_0}^\dagger}_F = \norm{\sum_{i=R_0+1}^R x_i x_i^\dagger}_F \le \sum_{i=R_0+1}^R \norm{x_i}_F^2 \le \alpha_\textrm{screen} \tau  , 
\end{equation}
where $X_{R_0} = [x_1, ..., x_{R_0}]$ is the first $R_0$ columns of $X$, and $\alpha_\textrm{screen}$ is a hyperparameter specifying how much of the error budget $\tau$ we're allowed to spend discarding columns of $X$. We typically set $\alpha_\textrm{screen} = 0.7$.

Performing this norm screening comes at minimal added cost because computing $\norm{x_i}_2^2 = \langle x_i, x_i \rangle$ is already a necessary step of the truncation algorithm. The only overhead is in ordering $\norm{x_i}_F$ to select which columns to discard.
Moreover, the upper bound $\sum_{i=R_0+1}^R \norm{x_i}_F^2$ from Eqn. \eqref{eqn:norm-screening} can be computed at no added cost, and with this quantity we can adaptively select the remaining truncation error $\tau_0$ allowed for the rest of the algorithm. In particular, the tolerance
\begin{equation}
    \tau_0 = \tau - \sum_{i=R_0+1}^R \norm{x_i}_F^2
\end{equation}
for the rest of the algorithm ensures the overall truncation error is below $\tau$.

\subsection{Partitioning error between SVD truncation and TT arithmetic}

One must decide how to partition the remaining error budget $\tau_0$ between the two sources of error: SVD truncation of the density matrix and inexact TT arithmetic. Writing $\tau_0 = \tau_\textrm{SVD} + \tau_\textrm{TT}$ to encode this partition, we typically allocate more of the budget towards SVD truncation as $\tau_\textrm{SVD} = 0.7 \tau_0$. The more allocated to the SVD, the fewer columns of $\tilde{X}$ need to be formed, each of which is a linear combination of $O(R_0)$ TT tensors. This is a heuristic choice, however, that results in potentially higher bond dimensions for the columns $\tilde{x}_i$. The optimal division of error will depend on the application.

\subsection{Linear combinations via randomized TT/MPS rounding}
\label{sec:linear-combos-randomized}

Once the number of columns $r$ of $\tilde{X}$ is determined, each $\tilde{x}_i = \sum_{j=1}^{R_0} V_r(i,j) x_j$ must be calculated in TT/MPS format.
One approach to perform each linear combination as a sequence of $R_0-1$ additions, with each addition followed by compression via TT-SVD. Refer to Supplement Material \ref{SM-alg:iterative-TT-SVD} for more details.
Although this approach allows for deterministic error control, it becomes very expensive when the number of terms and/or their bond dimensions are large.
In these cases, it is often much more efficient to use randomized methods \cite{Daas2021-Randomized-TT-Rounding, Daas2025-AdaptiveRT} to compute each $\tilde{x}_i$.

To compress the sum $\tilde{x}_i^\textrm{exact} = \sum_{j=1}^{R_0} V(j,i) x_j$ using the \texttt{Randomize} - \texttt{then} - \texttt{Orthogonalize} scheme of \cite{Daas2021-Randomized-TT-Rounding}, one first sketches each $x_j$, $j = 1, ..., R_0$, using another tensor train $\omega$ with stochastically generated cores.
This amounts to computing a sequence of \textit{partial contractions} $P_{j,k}$, $k = 2, ..., d$, c.f. \cite{Daas2021-Randomized-TT-Rounding} for details.
The compressed tensor $\tilde{x}_i$ is then computed from $\{P_{i,j}\}$ and the coefficients $V(:,i)$.

Randomized rounding techniques are particularly nice within the truncation routine because we must compute multiple linear combinations of the same set of vectors.
We use a single sketching tensor $\omega$ to perform all of these linear combinations rather than using a different random tensor per sum.
The associated partial contractions comprise a non-trivial portion of the randomized truncation algorithm, so re-using the sketching tensor can result in notable speed-up.

The bond dimensions $b^{(\omega)}$ of the sketching tensor $\omega$ define the maximal bond dimensions for the compressed representation of $\tilde{x}_i$.
Empirically, the bond dimensions of each $\tilde{x}_i$ change slowly between timesteps, so we often have a good approximation of what to select for $b^{(\omega)}$, say 1.2x the largest bond dimension of the components $x_1, ..., x_{R_0}$ in the sum.
Adaptive rank selection schemes are also available \cite{Daas2025-AdaptiveRT}, though not used here.

After randomization, an SVD truncation sweep is applied reduce the bond dimensions below $b^{(\omega)}$ if possible.
When performing this sweep, we set the SVD truncation tolerance per core to 
\begin{equation}
    \tau_\textrm{SVD}^{(i)} = \frac{1}{\sqrt{d-1}}\left(-\sigma_i + \sqrt{\sigma_i^2 + \frac{\tau_\textrm{TT}^2}{r^2}}\right) 
\end{equation}
because, disregarding error due to randomization, this ensures that 
\begin{equation*}
\norm{ \tilde{x}_i \tilde{x}_i^\dagger - \tilde{x}_i^\textrm{exact}(\tilde{x}_i^\textrm{exact})^\dagger }_F \le \frac{\tau_\textrm{TT}}{r}  . 
\end{equation*}
The overall truncation error is then bounded as $\lVert \tilde{X}\tilde{X}^\dagger - \tilde{X}^\textrm{exact}(\tilde{X}^\textrm{exact})^\dagger \rVert _F \le \tau_\textrm{TT}$.

\begin{algorithm}[t!]
\caption{Density matrix rank-compression scheme in TT/MPS format}
\label{alg:density-matrix-truncation}
\begin{algorithmic}[1]
\Require Factor matrix $X \in \C^{N \times R}$ with columns $x_i \in \C^N$ represented in TT/MPS format; \ Error tolerance $\tau$.
\vspace{0.35\baselineskip}
\Ensure Lower-rank factor matrix $\tilde{X} \in \C^{N \times r}$ with $\lVert XX^\dagger - \tilde{X}\tilde{X}^\dagger\rVert_F \le \tau$ and $r < R$
\vspace{0.35\baselineskip}
\Procedure{TT-Compress}{$X, \tau$}
    \vspace{0.5\baselineskip}
    \Statex $\quad$ \textbf{Norm Screening}
    \State Sort columns $x_i$ in decreasing order of norm
    \State $\tau_\textrm{screen} = \alpha_{\textrm{screen}} \ \tau$
    \Comment Allocate some error to norm screening
    \State $r_{\textrm{screen}} = \textrm{arg} \min \left\{ k \ \big| \ \sum_{i=k+1}^r \norm{x_i}^2  \le \tau_{\textrm{screen}} \right \}$
    \State $X \gets X(:,:r_{\textrm{screen}}) = [x_1, ..., x_{r_{\textrm{screen}}}]$
    \State $\tau \gets \tau - \sum_{i=r_\textrm{screen}+1}^R \norm{x_i}^2$
    \vspace{0.5\baselineskip}
    \Statex $\quad$ \textbf{Selection of Truncation Rank}
    \State $\tau_\textrm{SVD} = \alpha_\textrm{SVD} \tau$ \Comment{Allocate some error to SVD truncation}
    \State $W_{ij} = \langle x_i, x_j \rangle$ \ \textrm{for} $i,j = 1, ..., R$
    \State $W = V \Sigma^2 V^\dagger$ 
    \State $r = \textrm{arg} \min \left\{ k \ \big| \ \sum_{i=k+1}^R \sigma_i^2  \le \tau_\textrm{SVD}\right\}$
    \State $\tau \gets \tau - \sum_{i=r+1}^R \sigma_i^2$
    \vspace{0.5\baselineskip}
    \Statex $\quad$ \textbf{Calculation of $\tilde{X}$ in TT/MPS format using randomized rounding}
    \State $b^{(\omega)} = 1.2 \times (\textrm{maximum bond dimension out of all} \ x_1, ..., x_{R_0})$ 
    \State $\omega \gets$ sketching tensor with max bond dimension $b^{(\omega)}$
    \For {$i = 1, 2, ..., R_0$}
        \State $\{P_{j,k}\}_{k=2}^d = \texttt{Partial-Contraction-RL}(x_i, \omega)$
	\EndFor
    \For {$i = 1, 2, ..., r$}
        \State Compute $\tilde{x}_i$ from $V(:,i)$ and $\{P_{j,k}\}$
        \State $\tilde{x}_i \gets \texttt{SVD-Sweep-LR}(\tilde{x}_i, \ \tau_\textrm{TT}^{(i)})$ with $\tau_\textrm{TT}^{(i)} = -\sigma_i + \left(\sigma_i^2 + \tau^2 \ / \ r^2\right)^{1/2}$
	\EndFor
\EndProcedure
\end{algorithmic}
\end{algorithm}

\subsection{Pseudocode}

For clarity, Algorithm \ref{alg:density-matrix-truncation} provides the pseudocode for the rank truncation algorithm, incorporating the aspects discussed in this section.

%
%
\section{Numerical Experiments}
\label{sec:numerical-experiments}

We now present three numerical experiments to demonstrate the capabilities of the proposed scheme. 
Our implementation builds on the TT-Toolbox library \cite{TT-Toolbox} in \textsc{Matlab}.
All experiments were run using 16 cores of an AMD EPYC 7702 CPU (2 GHz) within the Advanced Research Computing (ARC) cluster at Virginia Polytechnic Institute and State University.

\subsection{Dissipative Spin Chain}

We first consider the spin-1/2 XX Heisenberg chain with dissipation at each site. This model describes a chain of $d$ magnetic spins whose individual Hilbert spaces are spanned by the states $\ket{\uparrow} = [1,0]^\top$ and $\ket{\downarrow} = [0,1]^\top$. The spins have nearest neighbors interactions modeled by the Hamiltonian
\begin{equation}
    H = \sum_{j=1}^{d-1} (\sigma_j^+ \sigma_{j+1}^- + \sigma_j^- \sigma_{j+1}^+) .
\end{equation}
Here, $\sigma^+ = \ket{\uparrow}\bra{\downarrow}$ is the single-site raising operator, $\sigma^- = (\sigma^+)^\dagger = \ket{\downarrow}\bra{\uparrow}$ is the associated lowering operator, and adding the subscript $\sigma_{j}^+$ for instance means the operator is applied to the $j$-th site: $\sigma_{j}^+ = I_{2\times 2}^{\otimes d-j} \otimes \sigma^+ \otimes I_{2\times 2}^{\otimes j-1} $.
Governed by the Lindblad equation, these spins are dissipative in the sense that the jump operators $L_j = \sigma_j^-/\sqrt{T}_\textrm{decay} $, $j = 1, ..., d$, drive the systems towards the pure system $\ket{\downarrow} \otimes \ket{\downarrow} \otimes \cdots \otimes \ket{\downarrow}$ over time.

\subsubsection{Convergence on Small Spin Chains} 

\begin{figure}[t!]
    \centering
    \includegraphics[width=\textwidth]{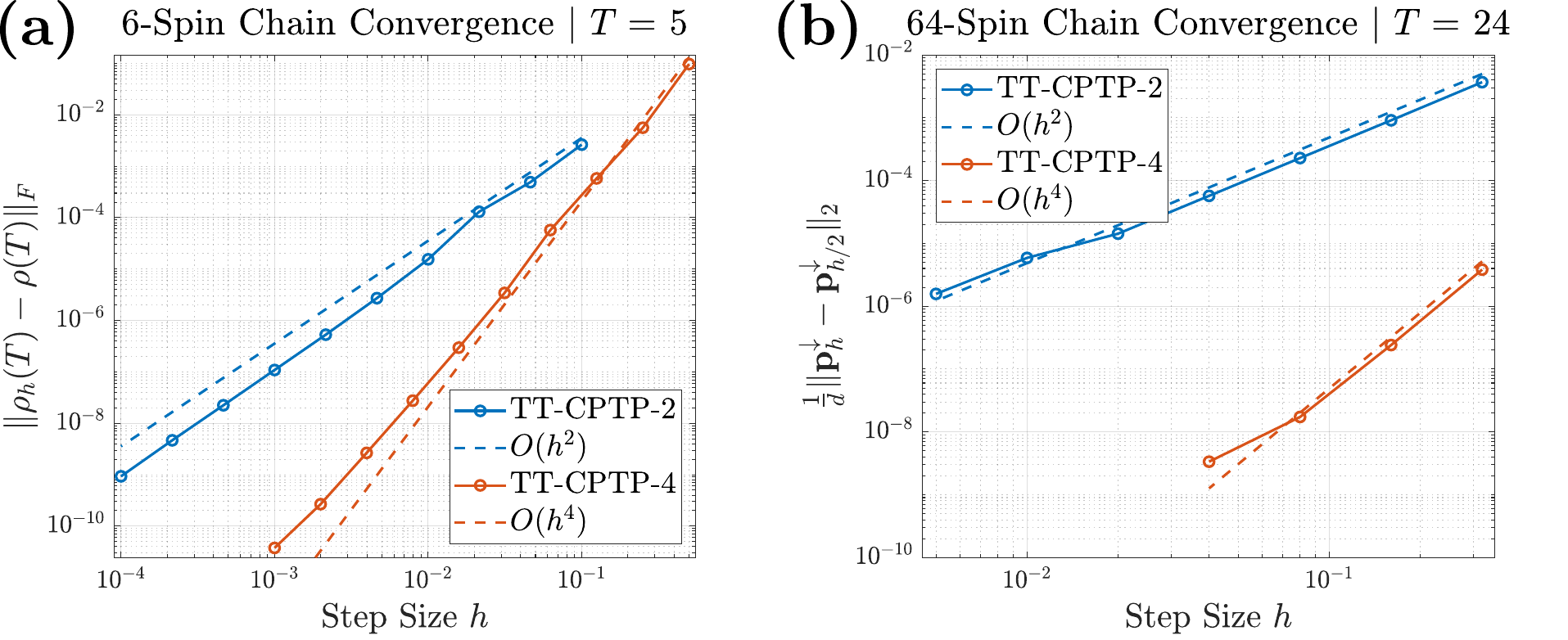}
    \caption{Convergence of the TT/MPS-based integrator applied to the dissipative 1D Heisenberg model. \textbf{(a)} shows the convergence in error $\norm{\rho_h(T) - \rho(T)}_F$ at $T = 5$ for $d = 6$ spins, for which a reference solution $\rho(T)$ can be computed using matrix-exponentiation of the Lindbladian. 
    Subfigure \textbf{(b)} shows the behavior on a larger system (64 spins) whose full Hilbert space has size $> 10^{19}$. Please refer to the text for more in-depth explanations.
    }
    \label{fig:heisenberg-chain-convergence} 
\end{figure}

As a first demonstration, we compute the error of the method on small systems ($d = 6$ spins) for which reference solutions can be obtained via matrix-exponentiation of the Lindbladian operator. We take the initial state to be $\rho(0) = \psi_0 \psi_0^\dagger$ with $\psi_0 = \ket{\uparrow \downarrow \downarrow \downarrow \downarrow \uparrow}$. The dissipation timescale is set at $T_\textrm{decay} = 20$, and we simulate the system for $T = 5$ units of times using our second and fourth order low-rank integrators using step sizes $h$ ranging from $10^{-4}$ to $5 \cdot 10^{-5}$. Operator splitting (e.g. TEBD) is employed to perform time-evolution with respect to the effective Hamiltonian $H + \frac{1}{2i} \sum_{j=1}^d L_j^\dagger L_j$, and randomized TT rounding is used within our density matrix compression scheme. We set the truncation error per timestep to $h^{p+1}$ where $p$ denotes the method order.

Figure \ref{fig:heisenberg-chain-convergence}\textbf{(a)} plots the convergence of the method in terms of $\norm{\rho_h(T) - \rho(T)}_F$, where $\rho_h(T)$ denotes the numerical solution at time $T$ obtained when using step size $h$, and $\rho(T)$ denotes the reference solution.
We see the expected convergence rates for both the order 2 and order 4 methods for errors above $10^{-8}$, though the order 4 method converges more slowly thereafter.
This result likely stems from our implementation of the fourth-order TEBD integrator, whose convergence as a solver for the Schrodinger equation stagnates around $10^{-10}$.

For our simulations of this length $6$ spin chain, the rank of the density matrix does not exceed 8, regardless of method order and step size, though the rank may grow is the simulation proceeded to later times $t > 5$.
We remark on the density matrix rank simply to indicate we are in a ``low-rank'' regime, as $8 < 2^6 = 64$, the size of the full Hilbert space of the 6 spins.

\subsubsection{Behavior for Larger Systems} 

When representing the density matrix as $\rho = VV^\dagger$ where $V$ is a dense matrix, the columns themselves have $2^d$ elements, which would be on the order of $10^{19}$ for 64 spins. Each of these dense vectors would need more than \textit{ten exabytes} to store in memory in single precision. 

Depending on the initial state, tensor train compression of the columns can make such simulations feasible using only modest compute resources. We demonstrate this result by simulating a system of $d = 64$ spins starting from pure state with all spins down aside from sites 8 and 48: $\psi_0 = \sigma_{8}^x \sigma_{48}^x \ket{\downarrow \downarrow \cdots \downarrow}$, where $\sigma_j^x$ flips the spin at site $j$. 
%
%
We set the dissipation timescale to $T_\textrm{decay} = 20$, and we simulate the system using both our second-order and fourth-order CPTP schemes with step sizes ranging from $h = 2.5 \cdot 10^{-3}$ to $4 \cdot 10^{-2}$. Operator splitting (TEBD) is again used to perform time-evolution w.r.t. the effective Hamiltonian. We set the truncation tolerances such that error $\tau = 10^{-5}/h$ and $\tau = 10^{-9}/h$ are introduced to the density matrix per timestep for the order-2 and order-4 methods, respectively. We used the smaller truncation threshold for the order-4 method because it converges so quickly.

A reference solution is not available at systems of this size, so we instead compare the output of the integrator when the step size is decreased by a half. For each step size $h$, we compute the probabilities $\mathbf{p}^\downarrow_h \in [0,1]^d$ at time $T = 20$ of measuring each spin to be in the state $\ket{\downarrow}$. This is the vector whose component $p_h(j)$ is a diagonal element of the reduced density matrix $\rho_h^{(j)}$ obtained by tracing out spins $k \neq j$:
\begin{equation}
    \mathbf{p}^{\downarrow}(j) := \bra{\downarrow}\rho^{(j)}_h\ket{\downarrow}, \quad \rho_h^{(j)} := \trace_{k \neq j}(\rho_h) .
\end{equation}
To demonstrate convergence of our methods, Figure \ref{fig:heisenberg-chain-convergence}\textbf{(b)} plots the deviation $\Delta_h := \lVert \mathbf{p}^\downarrow_h - \mathbf{p}^\downarrow_{h/2} \rVert$ in these probabilities when refining the step size by a factor of 2. These quantities $\Delta_h$ should scale as $O(h^{p})$ where $p$ is the order of the method, and this is indeed what is observed numerically, at least up to the truncation tolerance $\tau = 10^{-5}/h$ for the order-2 method and $\tau = 10^{-9}/h$ for the order-4 method. These results verify the methods converge with the right order. 

\begin{figure}[t!]
    \centering
    \includegraphics[width=\textwidth]{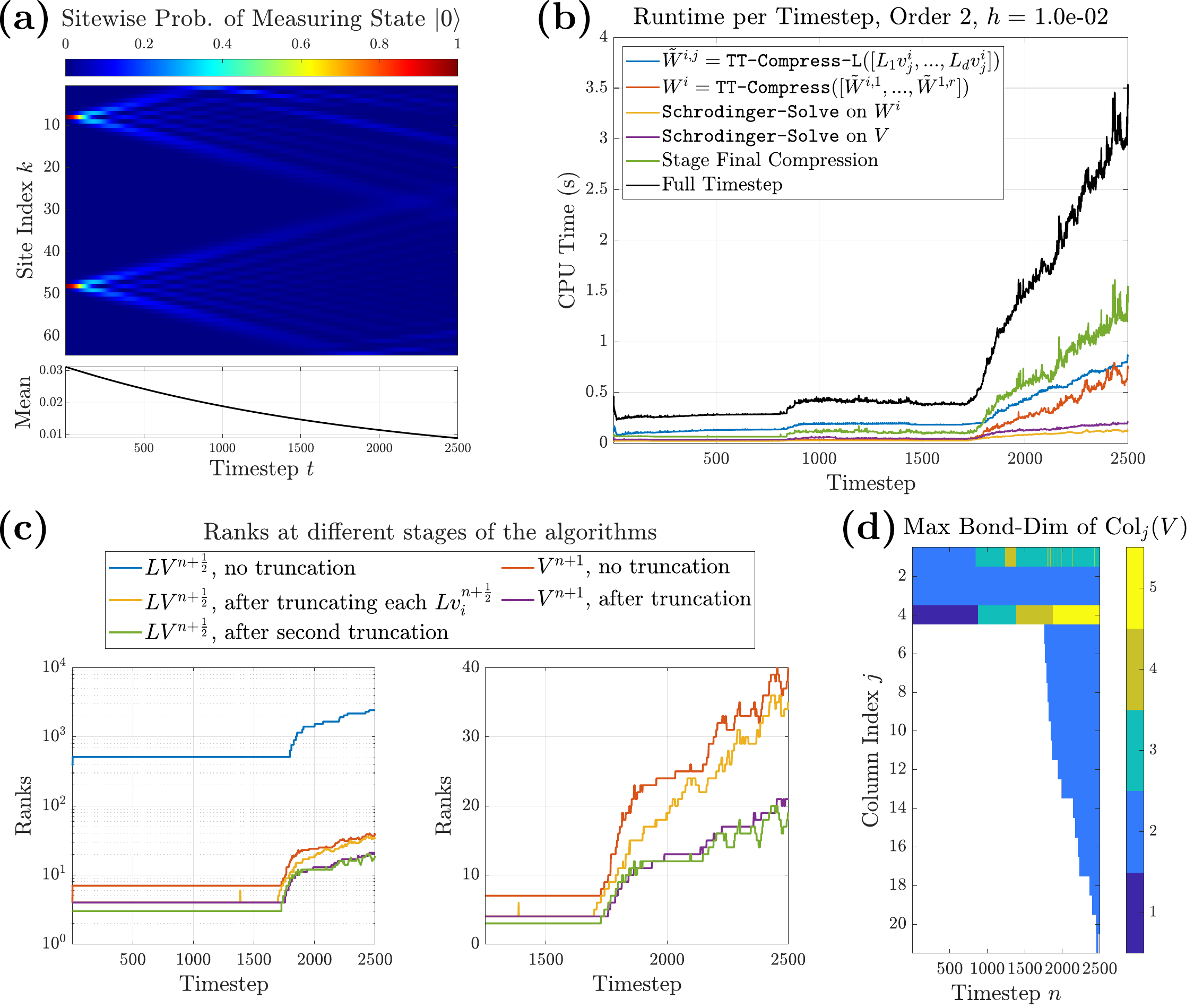}
    \caption{Low-rank dynamics of the dissipative 1D Heisenberg model with $d = 64$ spins, whose full Hilbert space has size $> 10^{19}$. The larger system starts in the state $\psi_0 = \sigma_{8}^x \sigma_{48}^x \ket{\downarrow \downarrow \cdots \downarrow}$, where $\sigma_j^x$ flips the spin at site $j$. The dissipative timescale is set to $T_\textrm{decay} = 20$.
    \textbf{(a)}-\textbf{(d)} record various run statistics when using our second-order CPTP scheme based on the midpoint rule with step size $h = 0.01$.
    Please refer to the text for more in-depth explanations of the subfigures.
    }
    \label{fig:heisenberg-chain-numerics} 
\end{figure}

Figure \ref{fig:heisenberg-chain-numerics}\textbf{(a)}-\textbf{(d)} aggregates various information recorded during the simulation using the method order $p = 2$ with step size $h = 0.01$. Subfigure \textbf{(a)} shows the probabilities $p_h^\downarrow(j)$ of measuring each spin in state $\ket{\downarrow}$, as well as its mean value $\frac{1}{d}\sum_{j=1}^d p^\downarrow_h(j)$ over the course of the simulation. The initial excitations (flip spins) at sites 8 and 48 propagate outwards as wave fronts from these sources, and they don't interact until around time $t = 20$.

Subfigure \textbf{(b)} plots the runtime of the simulation by timestep, broken down into different subroutines. Each timestep takes less than a second while the wavefronts interact, with the three types of truncations being the primary commutation cost. For clarity --- the blue curve denotes truncating matrices of the form $[L_1 v_j, ..., L_d v_j]$, with $v_i$ being a \textit{single column} of either $V^n$ or the middle stage $V^{n+1/2}$. 
After these matrices are truncated individually (c.f. Section \ref{sec:jump-operator-compression}), the remaining columns are passed to the generic truncation algorithm, whose runtime is plotted in red. This two-step process builds lower rank approximations $W^n$ and $W^{n+1/2}$ of $W_\textrm{full}^n := [L_1 V^n, ..., L_d V^n]$ and $W_\textrm{full}^{n+1/2} := [L_1 V^{n+1/2}, ..., L_d V^{n+1/2}]$, respectively, which are then used in forming $V^{n+1/2}$ and $V^{n+1}$. The green curves show the cumulative amount of time needed to truncate $V^{n+1/2}$ and $V^{n+1}$ after we've computed $W^n$ and $W^{n+1/2}$.

The runtime increases significantly once the wavefronts reach each other, with the truncations of $V^{n+1}$ and $W^{n+1/2}$ being the primary bottleneck despite the randomized TT rounding. Figure \ref{fig:heisenberg-chain-numerics}\textbf{(c)} shows the sizes (number of columns) of different factor matrices appearing during the truncation of $V^{(n+1)}$. Without truncation, the factors $W^{n+1/2}$ are large, with upwards of 1440 columns towards the end of the simulation. These matrices luckily have many small singular values that can be identified by truncating the matrices $[L_1 v_i^{(n+1/2)}, ..., L_d v^{(n+1/2)}_i]$ individually. These truncations, which take advantage of shared core structure when performing inner products and linear combinations (c.f. Section \ref{sec:jump-operator-compression}), reduce the number of columns by over an order of magnitude.
Truncation of the remaining columns via the generic truncation algorithm further cut the number of columns of $W^{n+\frac{1}{2}}$ in half.

Finally, Figure \ref{fig:heisenberg-chain-numerics}\textbf{(d)} shows the maximum bond dimension of each column of $V^n$ throughout the simulation. These bond dimensions remain small throughout the simulation, never exceeding 5, which amounts to huge compression of the columns of $V^n$. Indeed, at most $d n b^2 = 3200$ doubles are needed to represent a $d=64$ site MPS with max bond dimension $b \le 5$ with physical dimension $n = 2$, whereas $2^{64} > 10^{19}$ doubles would be needed to represent the a dense vector in the full Hilbert space.

\subsection{Mock Quantum Circuit Simulation}
\label{sec:mock-quantum-circuit-sim}

As a second example, we consider a collection of transmon qubits that weakly interact via the Jaynes-Cummings coupling. A sequence of two qubit SWAP gates is applied to these qubits via a piecewise constant control Hamiltonians, resulting in a mock simulation of a quantum circuit. It is a ``mock'' simulation in the sense that control Hamiltonian does not take the same form as those in actual quantum hardware, wherein controls are high frequency pulses optimized to implement quantum gates \cite{Petersson2021-QControl, Lee2025-HOHO}.

\subsubsection{System Description}

Following \cite{Petersson2021-QControl}, we consider a Hamiltonian describing a collection of transmon qubits that interact with each other via their coupling to resonator buses. When these buses are adiabatically eliminated and the qubits are moded as 2-level systems, the resulting Hamiltonian under the rotating wave approximation is
\begin{equation}
    H = \sum_{\langle p, q \rangle} J_{pq}(a_p a_q^\dagger + a_p^\dagger a_q) + H_\textrm{control}(t) .
\end{equation}
Here, $a_q = I^{\otimes (q-1)} \otimes \ket{0}\bra{1} \otimes I^{q-p+1}$ denotes the lowering operator for qubit $p$, and $J_{pq}$ is the Jaynes-Cummings coupling strength between qubits $p$ and $q$.
The sum is taken over pairs $\langle p, q \rangle$ determined by an underlying qubit layout, in our case a 25-qubit portion of a heavy-hex lattice. The layout, shown in Figure \ref{fig:heavy-hex-diagram}(a), is a two-dimension structure through which we snake the one-dimension MPS using the heuristically chosen qubit ordering in the figure. There are therefore long-range interactions encoded in the Hamiltonian, for instance between qubits 2 and 12, resulting in slightly higher MPO bond dimension than in the spin chain example.

The control Hamiltonian $H_\textrm{control}(t)$ is taken to be piecewise constant and implements a sequence of SWAP gates applied to different pairs of qubits in the system. It takes the form
\begin{equation}
\label{eqn:swap-gate-hamiltonian}    
     H_\textrm{control}(t) = \begin{cases}
        H_\textrm{SWAP}^{1,2} + H_\textrm{SWAP}^{10,11} + H_\textrm{SWAP}^{13,14} + H_\textrm{SWAP}^{23,24} 
        \ , & 
        t \in [0,T_\textrm{gate}) \\[0.5em]
        H_\textrm{SWAP}^{2,3} + H_\textrm{SWAP}^{9,10} + H_\textrm{SWAP}^{14,15} + H_\textrm{SWAP}^{22,23}
        \ , & 
        t \in [T_\textrm{gate},2T_\textrm{gate}) \\[0.5em]
        H_\textrm{SWAP}^{3,4} + H_\textrm{SWAP}^{8,9} + H_\textrm{SWAP}^{15,16} + H_\textrm{SWAP}^{21,22}
        \ , & 
        t \in [2T_\textrm{gate},3T_\textrm{gate})\\[0.5em]
        H_\textrm{SWAP}^{4,5} + H_\textrm{SWAP}^{7,8} + H_\textrm{SWAP}^{16,17} + H_\textrm{SWAP}^{20,21}
        \ , & 
        t \in [3T_\textrm{gate},4T_\textrm{gate}) \\[0.5em]
        H_\textrm{SWAP}^{5,6} + H_\textrm{SWAP}^{17,18} 
        \ , & 
        t \in [4T_\textrm{gate},5T_\textrm{gate})
    \end{cases}
\end{equation}
Here, $H_\textrm{SWAP}^{pq}$ is a time-independent Hamiltonian that, in the presence of the coupling $H_{JC}^{pq} = J_{pq}(a_p a_q^\dagger + a_p^\dagger a_q)$ between qubits $p$ and $q$, implements the unitary 
\begin{equation}
    U_\textrm{SWAP}^{pq} \ket{b_1 ... b_{p-1} b_p ... b_{q-1} b_q ... b_d} = \ket{b_1 ... b_{p-1} \textcolor{red}{b_q} ... b_{q-1} \textcolor{red}{b_p} ... b_d} \ , \quad b_i \in \{0,1\}
\end{equation}
which swaps the states of qubits $p$ and $q$. This is to say that
\begin{align}
    H_\textrm{SWAP}^{pq} &= -\frac{1}{i T_\textrm{gate}} \log(U_\textrm{SWAP}^{pq}) - H_{JC}^{pq}\ ,  
    \quad \textrm{so that}
    \\
     U_\textrm{SWAP}^{pq} &= \exp \left\{ -i T_\textrm{gate} (H_\textrm{SWAP}^{pq} + H_\textrm{JC}^{pq})\right\} \ . 
\end{align}
Essentially applying the Hamiltonian $H_\textrm{SWAP}^{pq}$ negates the coupling between the two active qubits and applies the swap gate to them.

The system is a modeled as an open quantum system governed by the Lindblad equation with two types of jump operators: \textit{decay} and \textit{dephasing}. These have the form
\begin{equation}
    \label{eqn:decay-and-dephasing}
    L_p^\textrm{decay} = \frac{1}{\sqrt{T_p^\textrm{decay}}} a_p , \quad \textrm{and} \quad  L_p^\textrm{dephase} = \frac{1}{\sqrt{T_p^\textrm{dephase}}} a_p^\dagger a_p,
\end{equation}
where the index $p$ indicates dependence on the qubit.

\subsubsection{Further Details}

\begin{figure}[t!]
    \centering
    \includegraphics[width=\textwidth]{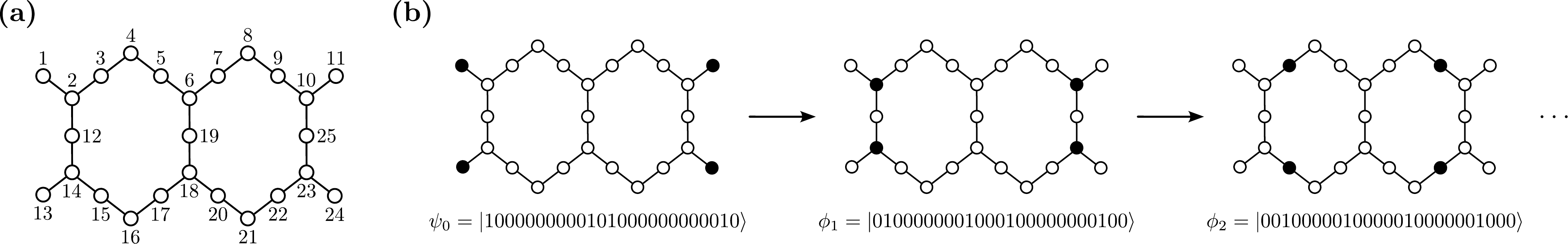}
    \caption{25-qubit portion of a larger heavy-hex lattice. Subfigure \textbf{(a)} shows our heuristic ordering imposed on the qubits when representing their state using an 25-site MPS. The Hilbert space has size $2^{25} \approx 3.4 \cdot 10^7$.
    Subfigure \textbf{(b)} shows the first three states $\phi_0$, $\phi_1$, and $\phi_2$ that the uncoupled system ($J_{pq} = 0$) would process through without due to the gate Hamiltonians. Open and filled-in nodes in diagram denote states $\ket{0}$ and $\ket{1}$, respectively, for the associated qubit.
    }
    \label{fig:heavy-hex-diagram}
\end{figure}

about our mock quantum circuit simulation.

\begin{table}[b!]
\centering
\begin{tabular}{ | l | c | c | c |}
\hline
       & $J_{pq}$ & $T_p^\textrm{decay}$ & $T_p^\textrm{dephase}$ \\ 
 \hline 
 Mean  & 2.3 MHz$ / 2\pi$ & 95 $\mu$s & 100 $\mu$s\\  
 StDev & 5.4 MHz$ / 2\pi$ & 5 $\mu$s & 10 $\mu$s \\
 \hline
\end{tabular}
\caption{Device parameters for the mock quantum circuit circuit simulation. }
\label{tab:mock-quantum-circuit-params}
\end{table}

\textbf{Flow operators}: Operator-splitting for the flow operators $U_h = \exp\{-i hH_\textrm{eff}\}$ is not as straightforward due to these long-range interactions. We instead elect to use the Taylor-series approach (c.f. Section \ref{sec:TT-time-evolution}) to building $U_h$, for which we find 10-12 terms in the series suffices for convergence. 
Building these such flow operators is moderately expensive, requiring maybe 15-20 seconds of compute time, but this is small relative to the overall runtime of the method.

\textbf{Initial state and gate sequence}:  The system starts in the product state $\psi_0 = \ket{b_1 b_2 ... b_d}$ with $b_i = 1$ for $i \in \{ 1, 11, 13, 24\}$ and $b_i = 0$ otherwise. 
We then apply a sequence of five sets of gates via the Hamiltonian in Eqn. (\ref{eqn:swap-gate-hamiltonian}). Without the Jaynes-Cummings coupling and Lindbladian treatment, this Hamiltonian would cause the four 1's in the initial state to be transferred perfectly from the qubits 1, 11, 13, and 24 to qubits 6, 7, 18, and 20, respectively, at time $t = 5 T_\textrm{gate}$. The circuit would go through a sequence of intermediate states
\begin{align}
\begin{split}
    \phi_0 
    \ = \ \sigma_1^x \sigma_{11}^x \sigma_{13}^x \sigma_{24}^x \ket{0}^{\otimes 25} 
    \ &= \ \ket{1000000000101000000000010} = \psi_0,
    \\
    \phi_1 
    \ = \ \sigma_2^x \sigma_{10}^x \sigma_{14}^x \sigma_{23}^x \ket{0}^{\otimes 25} 
    \ &= \ \ket{0100000001000100000000100},
    \\
    \phi_2
    \ = \ \sigma_3^x \sigma_{9}^x \sigma_{15}^x \sigma_{22}^x \ket{0}^{\otimes 25} \ &= \ \ket{0010000010000010000001000},
    \\
    \textrm{etc}.
\end{split}
\end{align}
The first three states $\phi_0, \phi_1$, and $\phi_2$ are depicted in Figure \ref{fig:heavy-hex-diagram}(b) for reference. Due to the coupling $J_{pq} \neq 0$, the system never to reach these states exactly. We measure the extent to which this deviation occurs by computing the state populations
\begin{equation}
    \label{eqn:heavy-hex-phi-overlaps}
    \alpha_i(t) := \bra{\phi_i} \rho(t)\ket{\phi_i} = \sum_{j=1}^{r(t)} \big|\braket{v_j(t)}{\phi_i} \big|^2
\end{equation}
throughout the simulation.

\textbf{System parameters}: The gate time is fixed at $T_\textrm{gate} = 100$ nanoseconds. Qubit coupling strengthens, decay times, and dephasing times are generated i.i.d. from normal distributions with means and standard deviation as given in Table \ref{tab:mock-quantum-circuit-params}. We use the second order CPTP scheme based on the midpoint rule using a step size $h = 0.2$ nanoseconds and allowing $\tau = 10^{-5}$ truncation error into the density matrix per nanosecond (e.g. $10^{-6}$ error per timestep). Randomized rounding of MPO-MPS products and linear combinations is employed whenever the bond dimension of the uncompressed MPS exceeds 32.

\subsubsection{Numerical Results}

\begin{figure}[t!]
    \centering
    \includegraphics[width=\textwidth]{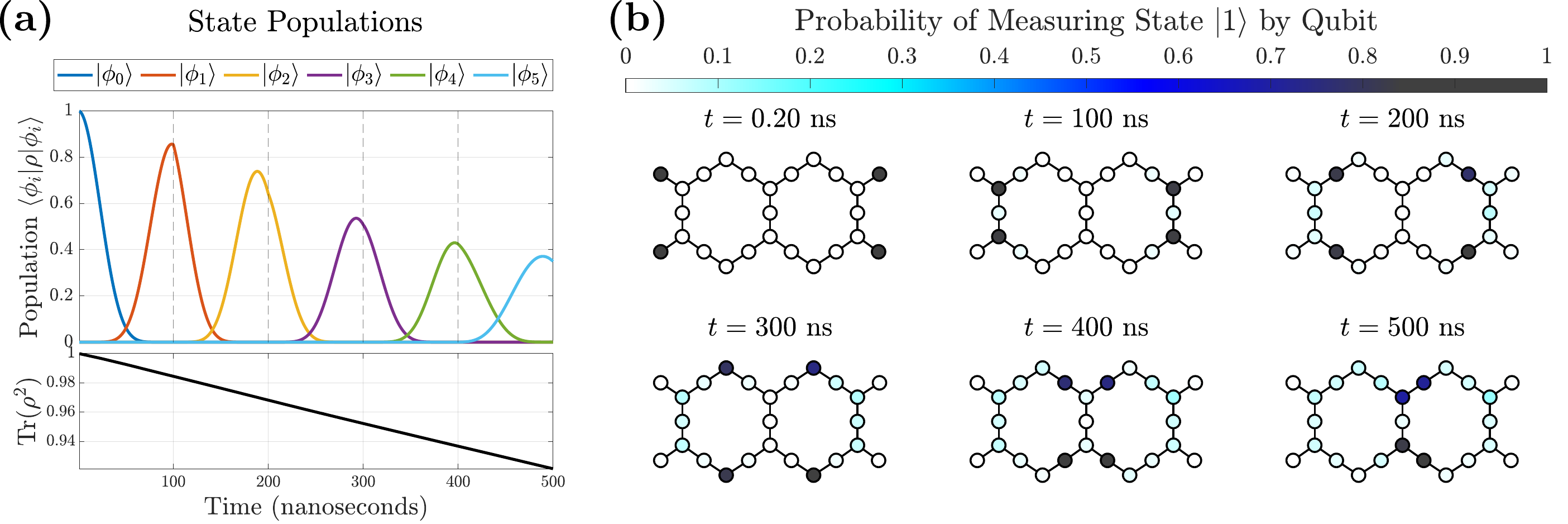}
    \caption{Dynamics of the 25-qubit heavy-hex lattice mock quantum circuit simulation.  
    Subfigure \textbf{(a)} plots the populations $\alpha_i(t) : =\bra{\phi_i} \rho(t) \ket{\phi_i}$ of the system in the sequence of states $\ket{\phi_0}, \ket{\phi_1}, ..., \ket{\phi_5}$ that the circuit would progress through without the Jaynes-Cummings coupling and qubit decay/dephasing.
    Each gate's control signal lasts 100 nanoseconds as indicated by the vertical dotted lines. The maximum populations in the target states achieved by the circuit decreases with each gate. 
    These errors are predominately due to the JC coupling, though the decay/dephasing do decrease the system's purity $\textrm{Tr}(\rho^2)$ to around 0.95 by the end of the simulation.
    \textbf{(b)} shows snapshots of the system state at times $t = h, \  T_\textrm{gate}, \  2 T_\textrm{gate}, ..., \ 5 T_\textrm{gate}$ overlaid on the heaxy-hex lattice.
    We show the probability of finding each qudit to be in its excited state $\ket{1}$ at these points in time.
    Only four qudits have excitations initially, and 
    these excitations propagate imperfectly through the system due to the Jaynes-Cummings coupling on top of the SWAP Hamiltonian sequence. 
    %
    }
    \label{fig:heavy-hex-dynamics}
\end{figure}

Figure \ref{fig:heavy-hex-dynamics} conveys data on the dynamics of the qudits during the mock circuit simulation.
Subfigure \textbf{(a)} plots the overlaps $\alpha_i(t)$ (c.f. Eqn. \ref{eqn:heavy-hex-phi-overlaps}) between the system state $\rho(t)$ and the sequence of states $\ket{\phi_i}$ the system would process through without the Jaynes-Cummings coupling and qubit decay/dephasing.
Without these sources of error, the circuit would satisfy $\alpha_j(i \cdot T_\textrm{gate}) = \delta_{ij}$, namely the system would be precisely in the pure state $\ketbra{\phi_i}{\phi_i}$ at time $t = i \cdot T_\textrm{gate}$.
Due to the coupling and Lindbladian treatment, this behavior is observed in our simulation, where one sees the maximal overlap between the system and the states $\ket{\phi_i}$ decreases with each gate.

Fig \ref{fig:heavy-hex-dynamics}\textbf{(b)} elucidates this behavior further by depicting the probability of measuring each qubit in its excited state $\ket{1}$ at times $t = \Delta, \ T_\textrm{gate}, \ 2 T_\textrm{gate}, ... , \ 5 T_\textrm{gate}$. These probabilities, denoted $p^{\ket{1}}(j)$ for $j = 1, ..., d$, are the (1,1) element of the reduced density matrices obtained by tracing out all but one qubit from the system, namely
\begin{equation}
    p^{\ket{1}}(j)  := \bra{1} \rho_j(t)\ket{1} \ , \enskip \rho_j(t) := \textrm{Tr}_{k \neq j} \big(\rho(t)\big).
\end{equation}
The system starts in the pure product state with its four outermost qubits in state $\ket{1}$ and the rest in state $\ket{0}$. 
Due to the sequence of SWAP gates, these ``excitations'' propagate between qubits over the course of the simulation, though this transmission is inexact due to the JC coupling and decay/dephasing.
As the excitations move through the circuit, some of their probability density is lost to nearby qubits that aren't involved in the gates.
By the time the final gate finishes, nearly all qubits have non-zero probability of being in their excited state. 
Note there is asymmetry in the excitation propagation due to the stochastically generated device parameters, which differ from qubit to qubit.

\begin{figure}[t!]
    \centering
    \includegraphics[width=\textwidth]{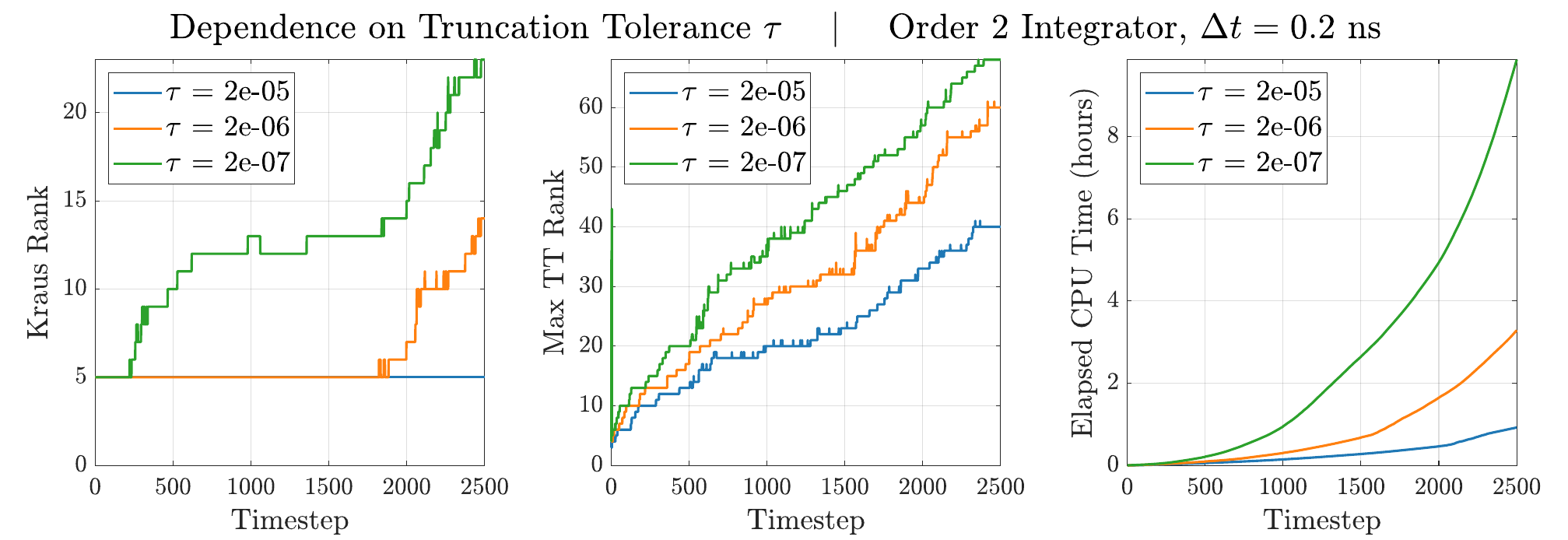}
    \caption{Run statistics of the heavy-hex lattice mock quantum circuit simulation: density matrix rank, maximum MPS bond dimension, and elapsed runtime by timestep when using three different rounding tolerances.
    }
    \label{fig:heavy-hex-run-stats}
\end{figure}

Figure \ref{fig:heavy-hex-run-stats} shows the dependence of the density matrix rank, maximum MPS bond dimension, and simulation runtime on the rounding tolerance $\tau$. 
All three metric increase as the tolerance decreases, with the smallest tested tolerance $\tau = 2 \cdot 10^{-7}$ using far more computational resources than the rest.
At this smallest tolerance, one must resolve upwards of 25 columns of the factor matrix $V$, at higher bond dimensions, resulting in this simulation lasting 3x longer than the simulation at the middle tolerance $\tau = 2 \cdot 10^{-6}$.

\subsection{Qudit-Resonator Chain}

As a final example, we again consider a systems of transmon qudits, but this time we also model resonator buses through which qubit interactions are mediated. The system consists of $d$ subsystems with physical dimensions $n_1, n_2, ..., n_d$, with each subsystem being either a qudit or a resonator. The overall Hamiltonian $H$ is the sum of the time-independent \textit{device Hamiltonian} $H_d$ and the time-dependent \textit{control Hamiltonian} $H_c(t)$.

Viewed in the rotating frame, the device Hamiltonian is given by 
\begin{equation}
    \label{eqn:circuit-qed-device-hamiltonian}
    H_d = \sum_{q} \left((\omega_q - \omega_q^\textrm{rot}) a_q^\dagger a_q - \frac{\xi_q}{2}a_q^\dagger a_q^\dagger a_q a_q\right) - \sum_{\langle p, q \rangle} \xi_{pq} a_p^\dagger a_p a_q^\dagger a_q .
\end{equation}
As before, $a_q$ denotes the lowering operator for subsystem $q$
\begin{align}
    a_q &= I_{n_d} \otimes \cdots \otimes I_{n_{q+1}} \otimes A_q \otimes I_{n_{q-1}} \otimes \cdots \otimes I_{n_1} \ , 
    \enskip \textrm{with} \\
     A_q &= \sum_{j=1}^{n_q-1} \sqrt{j} \ \ketbra{j}{j+1} \ \in \ \R^{n_q \times n_q} .
\end{align}
$\omega_q$ is the frequency of subsystem $q$, $\omega_q^\textrm{rot}$ is the frequency of the rotating frame for subsystem $q$, $\xi_q$ is the self-Kerr coefficient of subsystem $q$, and $\xi_{pq}$ is the cross-Kerr coupling strength between subsystems $p$ and $q$.
It is standard to take $\omega_q^\textrm{rot} = \omega_q$ so that the terms $a_q^\dagger a_q$ vanish from the Hamiltonian.
The sum introducing the coupling is taken over pairs $\langle p, q \rangle$ determined by an underlying device layout, this time a chain of alternating qudits and resonators (c.f. Figure \ref{fig:qudit-resonator-chain-diagram}). There are both nearest neighbor-coupling between qudit and resonators, as well as two-step coupling between qudits. 

For our example, the device has 6 qudits and 5 resonators. Qudit subsystems have $n_\textrm{qudit} = 4$ levels modeling two states $\ket{0}$ and $\ket{1}$ used as a computational basis and two guard states $\ket{3}$ and $\ket{4}$. Resonators have a higher number of states, $n_\textrm{res} = 10$, making the overall Hilbert space have dimension $N = n_1 \cdots n_d = 4^6 \cdot 10^5 \approx 4 \cdot 10^{8}$.

Whereas the device Hamiltonian $H_d$ has non-zero elements only along its diagonal,  the control Hamiltonian $H_c$ introduces off-diagonal elements. In the rotating frame,
\begin{equation}
    H_c(t) = \sum_{q=1}^d \big( d_q(t) a_q + \bar{d}_q(t) a_q^\dagger\big).
\end{equation}
Each $d_q(t)$ is a complex-valued function denoting the control signal to transmon $q$. These controls are designed to implement quantum gates, see for instance \cite{Lee2025-HOHO}.

As in the previous example, the qudits experience both \textit{decay} and \textit{dephasing} jump operators (c.f. Eqn. \ref{eqn:decay-and-dephasing}). 
Resonators only experience decay, with their the decay timescale being much smaller than that of the qudits.
$T^\textrm{decay}_q$ and $T^\textrm{dephase}_q$ are drawn from normal distributions with means given in the table below and standard deviation being $1\%$ of the mean.
\begin{center}
\begin{tabular}{ | l | c | c | c | c |}
\hline
    & $T_p^\textrm{decay}$ ($\mu$s) & $T_p^\textrm{dephase}$ ($\mu$s) \\ 
 \hline 
 Qudits     & 95  & 50 \\  
 Resonators & 0.4  & $\infty$ \\
 \hline
\end{tabular}
\end{center}

\subsubsection{Control Signals}

\begin{figure}[t!]
    \centering
    \includegraphics[width=0.9\textwidth]{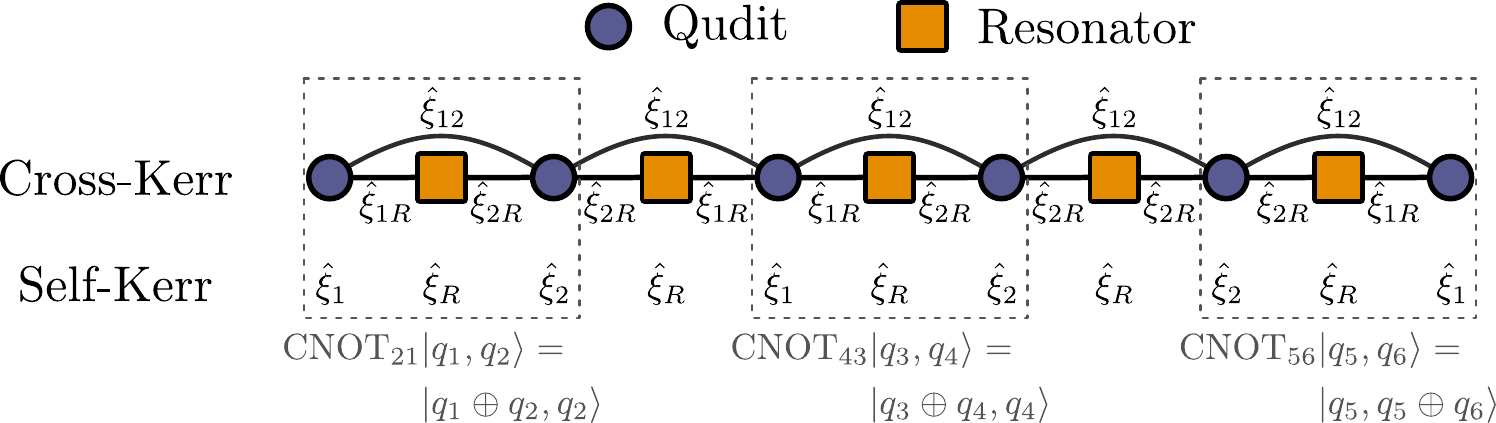}
    \caption{
    Chain of alternating qudits (4-level systems) and resonators (10-level systems), whose full Hilbert space has size $4^{6} \cdot 10^5 \approx 4 \cdot 10^8$.
    Cross-Kerr interactions (parameter $\xi_{pq}$ in Eqn \ref{eqn:circuit-qed-device-hamiltonian}) are denoted by the edges in the graph, with qudit-resonator interactions being stronger than qudit-qudit interactions.
    Parameter values $\hat{\xi}_q$, $\hat{\xi}_{pq}$, and $\hat{\xi}_{qR}$, as well as control signals $\hat{d}_q(t)$, are replicated from Example 4.2 of \cite{Lee2025-HOHO}, who optimize control signals to implement a controlled-not (CNOT) gate for a closed system with two qudits and one resonator.
    Each qudit in our larger 6-qudit chain is assigned either ``Type 1'' or ``Type 2'' to indicate it is a copy of either qudit 1 or qudit 2 from \cite{Lee2025-HOHO}. 
    All resonators use the same self-Kerr $\hat{\xi}_R$ values, but only resonators 1, 3, and 5 receive non-zero control signals $\hat{d}_R(t)$. 
    This choice results in approximate CNOT gates being applied to qudits 2 \& 1, qudits 4 \& 3, and qudits 5 \& 6 in our larger system. Note we write $\textrm{CNOT}_{21}$ in the diagram, e.g. with the 2 before the 1, to indicate the 2nd qudit is the control qudit of its associated CNOT gate.
    These gates are approximate due to the cross-Kerr interactions and our Lindbladian treatment of the system.
    }
    \label{fig:qudit-resonator-chain-diagram}
\end{figure}

In one of their examples, \cite{Lee2025-HOHO} optimize control signals for a three-transmon system (2x qudits and 1x resonator) to implement a controlled-not (CNOT) gate
\begin{equation}
    \textrm{CNOT} = \ketbra{000}{000} + \ketbra{001}{001} + \ketbra{100}{101} + \ketbra{101}{100},
\end{equation}
where $\ket{j_1j_Rj_2}$ indicates the product state with qudit $i$ in state $\ket{j_i}$ and the resonator in state $\ket{j_R}$.
They operator in the Schrodinger (closed) picture with device parameters given in the Table \ref{tab:qudit-resonator-circuit-params}. We've used $\hat{\cdot}$ notation to indicate these are parameters of the 3-transmon system of \cite{Lee2025-HOHO}. Similarly, let $\hat{d}_1, \hat{d}_R, \hat{d}_2$ denote their optimized control signals achieving the CNOT gate, in the sense that the solutions to the initial value problem $\frac{d}{dt}\ket{\psi} = -i(\hat{H}_d + \hat{H}_c(t)) \ket{\psi}, \quad \ket{\psi(0)} = \ket{\psi_0}$ satisfy $\ket{\psi(T_\textrm{gate})} = \textrm{CNOT} \ \ket{\psi_0}$. Their gate duration $T_\textrm{gate}$ is set to $550$ nanoseconds.

\begin{table}[t!]
\centering
\begin{tabular}{ | c | c | c | c | c | c | c |}
\hline
& $\hat{\xi}_1$ & $\hat{\xi}_2$ & $\hat{\xi}_R$ & $\hat{\xi}_{12}$ & $\hat{\xi}_{1R}$ & $\hat{\xi}_{2R}$ \\
\hline 
Value/2$\pi$ (GHz) & 0.220 & 0.225 & $2.83\cdot 10^{-3}$ & $10^{-6}$ & $2.49\cdot 10^{-3}$ & $2.52\cdot 10^{-3}$ \\
\hline
\end{tabular}
\caption{Device parameters for the qudit-resonator circuit simulation. }
\label{tab:qudit-resonator-circuit-params}
\end{table}

For our experiment, we duplicate the device parameters and optimized control signals from \cite{Lee2025-HOHO} three times to model three CNOT gates being performed simultaneously within our a larger system of 6 qudits and 5 resonators. 
These CNOT gates are applied between qudits 1 \& 2, 3 \& 4, and 5 \& 6, respectively.
Figure \ref{fig:qudit-resonator-chain-diagram} shows the assignment of parameters and control signals in the larger circuit. For instance, three qudits in our circuit have device parameters $\hat{\xi}_1$ and $\hat{\xi}_{1R}$, and they receive the control signal $\hat{d}_1(t)$.
Note also that the 2nd and 4th resonators receive no input signals, as CNOT gates are not appled between qudits 2 \& 3 and 4 \& 5, respectively.

The optimized control signals $\hat{d}_q$ from \cite{Lee2025-HOHO} are differentiable functions represented view B-splines. Our integrator is designed from time-independent, so we approximate the optimized signals $\hat{d}_{q_p}$ as piecewise constant functions, with value changing every $T_\textrm{signal} = 1$ nanosecond. This is to say that if transmon $p$ is assigned the optimized control signal $\hat{d}_q$, it instead receives the piecewise constant approximation
\begin{equation}
    d_p(t) = (\hat{d}_{q_p} \circ \mathfrak{T})(t) \ , \quad \textrm{with} \quad \mathfrak{T}(t) = \left \lfloor \frac{t}{T_\textrm{signal}} \right \rfloor \cdot T_\textrm{signal} + \frac{T_\textrm{signal}}{2}.
\end{equation}

\subsubsection{Flow Operators}

The Hamiltonian remains constant with each interval $\big[(k-1) T_\textrm{signal}, kT_\textrm{signal} \big)$ for $k = 1, ..., T_\textrm{gate}/T_\textrm{signal}$. 
Letting $H^{(k)}$ denote the Hamiltonian during the $k$-th time interval, we build an approximation of the flow operator $U^{(k)} := \exp\{-i h H^{(k)}\}$ via the operator splitting technique (TEBD).
The qudit-qudit cross-Kerr terms introduce interactions between MPS sites a distance 2 apart, so a bit more care is needed for this construction as compared to the spin-chain example.
Refer to supplemental material \ref{supp:qudit-resonator-tebd} for details.

\subsubsection{Further Experiment Details}

The system starts in a product state $    \ket{\psi_0} = \ket{00101000101}$. This is to say all of the resonators (the even numbered indices) start in their state $\ket{0}$, whereas the qudits start in either state $\ket{0}$ or $\ket{1}$. The initial qudits states were chosen to demonstrate the different actions of the CNOT gate. We use our second-order CPTP scheme using a timestep $h = 0.1$ nanoseconds, with the truncation tolerance set to $\tau = 5 \cdot 10^{-6}$ per timestep. As with the previous example, randomized rounding is used to compress MPO-MPS products as well as linear combinations of MPS.

\subsubsection{Numerical Results}

\begin{figure}[t!]
    \centering
    \includegraphics[width=\textwidth]{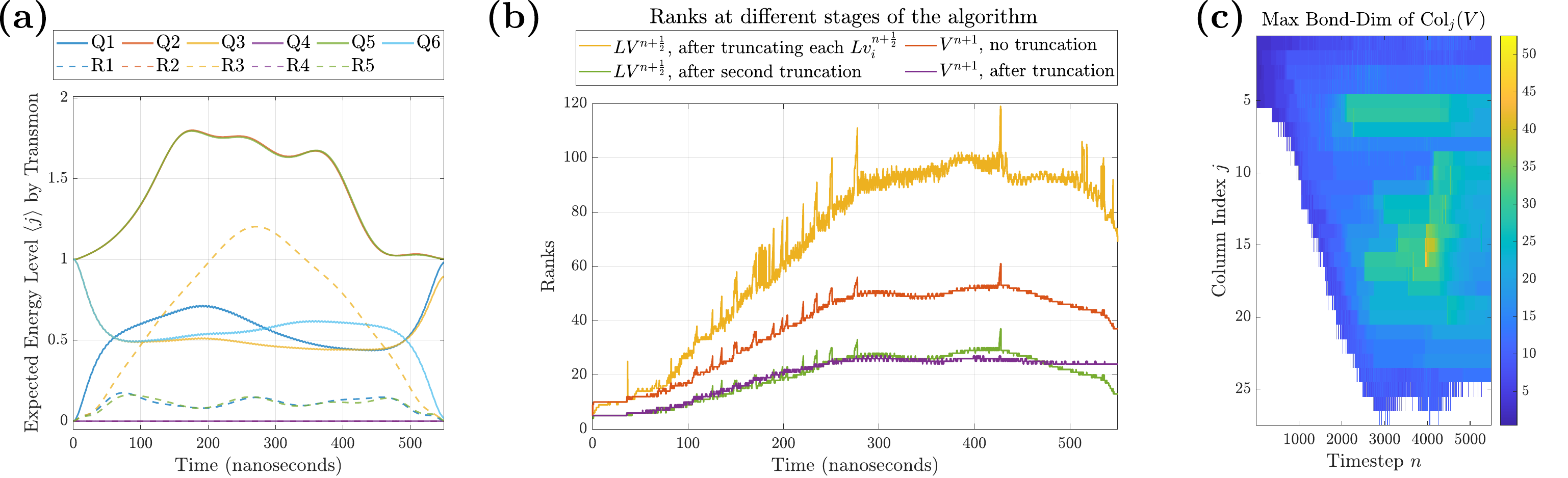}
    \caption{Low-rank dynamics of six qudits (4-level systems) and five resonators (10-level systems) aranged in a chain, as in Figure \ref{fig:qudit-resonator-chain-diagram}. Piecewise-constant control signals implement CNOT gates between qudits 1 \& 2, 3 \& 4, and 6 \& 5 over the course of 550 nanonseconds.
    We use the second-order CPTP integrator with timestep $h = 0.1$ nanosecond and $\tau = 5 \cdot 10^{-6}$ truncation error per timestep.
    Subfigure \textbf{(a)} shows the expected energy level $\langle j \rangle$ of each qudit and resonator.
    \textbf{(b)} plots the density matrix rank at different stages of the time-stepping algorithm at each point in time, and \textbf{(c)} indicates the maximum bond dimensions of each column $v_j(t)$ of the factor matrix $V(t)$. 
    Please refer to the text for more in-depth explanations of the subfigures.
    }
    \label{fig:qudit-resonator-numerics}
\end{figure}

Figure \ref{fig:qudit-resonator-numerics}\textbf{(a)} plots the evolution of the qudit resonator chain during the course of the simulation.
We specifically show the \textit{expected energy-level} for each transmon,
\begin{equation}
    \langle j_q \rangle (t) =  \sum_{j_1,j_2,...,j_d} j_q \bra{j_1 j_2 \cdots j_d} \rho(t) \ket{j_1 j_2 \cdots j_d} .
\end{equation}
The system starts in the state $\rho_0 = \ketbra{\psi_0}{\psi_0}$ with $\ket{\psi_0} = \ket{00101000101}$, so at $t = 0$, $  \langle j_1\rangle = 0$, $\langle j_2\rangle = 0$, $\langle j_3\rangle = 1$, $\langle j_4\rangle = 0$, etc.
If the control signals where optimized to account for transmon decay, dephasing and cross-talk, the final state should be the pure state $\hat{\ket{\psi}}\hat{\bra{\psi}}$ with
\begin{equation}
    \hat{\ket{\psi}} = \textrm{CNOT}_{21}\textrm{CNOT}_{43}\textrm{CNOT}_{56} \ket{\psi_0} = \ket{10101000100}.
\end{equation}
In that case, we'd read $\langle j_q \rangle \in \{0,1\}$ for all the qudits and $\langle j_q \rangle = 0$ for all of the resonators at time $t = T_\textrm{gate} = 500$ ns.
These controls were however optimized for a closed system that didn't account for all of the cross-talk within this larger quantum circuit, so the system does not terminate in the desire state.
Instead, from Figure \ref{fig:qudit-resonator-numerics}\textbf{(a)} we can see many of the qudits and resonators have expected energy level slightly off from $0$ or $1$ at the end of the simulation.

At intermediate times $t < T_\textrm{gate}$, the transmons exhibit non-trivial dynamics due to the control signals.
Of note, the resonators (dotted lines) generally remain at low energy levels, with only resonator 3 exceeding an expected energy level of 1. This limited activity of the resonators is likely what enables the system's state to remain low-rank, as indicated in Figure \ref{fig:qudit-resonator-numerics}\textbf{(b)}. Here we see the ranks of the density matrix at different stages of the time-integration algorithm, in the same style of the examples from the previous sections. The purple curve denotes the rank of the truncated density matrix at each point in time, which remains below 30 throughout the simulation. 
The ranks pre-truncation are higher, particularly for $\mathcal{L}_L V$, whose rank is upwards of 100 even after our initial truncation based on each column of $V$ individually.

Figure \ref{fig:qudit-resonator-numerics}\textbf{(c)} shows the maximum bond dimensions of the columns $v_j(t)$ of the factor matrix $V(t) = [v_1(t), ..., v_{r(t)}(t)]$ of the density matrix $\rho(t) = V(t) V(t)^\dagger$. Each row in the heatmap delineates a column $v_j(t)$ of $V(t)$. The color white at cell $(j,n)$ in the heatmap means that $V(t_n)$ had $r(t) < j$ columns.
The maximal bond dimensions generally remain low, with mean value never exceeding 25.
Some of the columns $v_j(t)$ do reach max bond dimension 52, which might still be considered ``low-rank''. Indeed, the MPS representations of these vectors need around 125,000 complex numbers, whereas their full representation as dense vectors would have more than $10^8$ elements. The MPS format therefore offers more than 3000x memory compression for these highest bond dimension vectors.
For the average case of bond dimension below 25, the compression is much higher, on the order of 10000x.

%
%
\section{Conclusion}

In this work we developed a family of low rank CPTP schemes for solving the Lindblad equation.
These schemes are exhibit two levels of low-rankedness, first by factorizing the density matrix as $\rho = VV^\dagger$, and second by representing the columns $v_i$ of $V$ in TT/MPS format.
This representation fits naturally into our recently developed Kraus-is-King scheme \cite{Appelo2024-kraus-is-king}, which boils down to a series of arithmetic operations on the columns of $V$.
As the primary bottleneck is the low-rank truncation of $V$ itself, we discuss in detail how to perform this compression efficiently in the TT/MPS format.
We demonstrate the capabilities of our methods on two representative open quantum systems from quantum computing and one open system from condensed matter physics, with Hilbert space dimensions greater than $10^8$ and $10^{19}$, respectively.

A natural extension of this work is to use our low-rank format within Lindblad schemes based on nested Picard iteration for time-dependent Hamiltonians \cite{Hu2025-Time-Dependent-Low-Rank-Lindblad}. The compression techniques developed in this paper should transfer directly to these integrators, which use the same arithmetic operators as the Kraus-is-King scheme.

On the implementation side, our methods can likely benefit from more careful management of parallel compute resources. Our implementation currently dedicates \textit{all available threads} to a single inner product or TT sum calculation, each of which involves operations on small matrices (due to low bound dimensions). With these matrices being so small, parallelism over inner products and linear combinations may be more efficient than BLAS/LAPACK parallelism within each of these routines.

%
%
\section*{Acknowledgements}

DA is supported by the U.S. Department of Energy, Office of Science, Advanced Scientific Computing Research (ASCR), under Award Number DE-SC0025424. This material is based upon work supported, in part, by the National Science Foundation under Grant No DMS-2436319 and Virginia Tech. 
YC is supported by DOE grant DE-SC0023164, AFOSR grant FA9550-25-1-0154 and Virginia Tech. This material is based upon work supported by the National Science Foundation under Grant No. DMS-2424139 while the second and third authors were in residence at the Simons Laufer Mathematical Sciences Institute in Berkeley, California, during the Fall 2025 semester.

%
%
\bibliographystyle{siamplain}
\bibliography{references,ref_lr}

%
%
\appendix
%
%
\section{Linear Combinations via TT-SVD}
\label{SM-sec:linear-combos-via-TT-SVD}

The primary cost of the density matrix compression scheme $X \in \C^{N \times R} \to \tilde{X} \in \C^{N \times r}$ is linear combinations performed in the TT/MPS format.
Following the notation of Section \ref{sec:accelerating-rounding}, we are interested in computing each column $\tilde{x}_i = \sum_{j=1}^{R_0}  V(j,i) x_j$, where the coefficients $V(j,i)$ have been determined via SVD truncation at the level of the density matrix.
One approach to performing each linear combination is via a sequence of $R_0 - 1$ sums in the TT format. For instance,
\begin{equation}
    \tilde{x}_i^{(j)} = \texttt{TT-SVD}\left(\tilde{x}_i^{(j-1)} + V(j,i) x_j, \ \tau_{ij}\right) \ , \enskip i = 1, 2, ..., r, 
\end{equation}
with
\begin{equation}
    \tilde{x}_i^{(1)} = V(1,i) x_1 \quad \textrm{and} \quad \tilde{x}_i := \tilde{x}_i^{(R_0-1)} \ . 
\end{equation}
Truncation makes each sum inexact, introducing perturbations
\begin{equation}
    \tilde{x}_i^{(j)} = \left[\tilde{x}_i^{(j-1)} + V(j,i) x_j\right] + e_{ij} \quad \textrm{with} \quad \norm{e_{ij}}_2 \le \tau_{ij}.
\end{equation}
The perturbation directions $e_i := \tilde{x}_i - \tilde{x}_i^\textrm{exact}$ in computing each $\tilde{x}_i$ therefore satisfies
\begin{equation}
    \norm{e_i}_2 \le \sum_{j=1}^{R_0} \norm{e_{ij}}_2 \le \sum_{j=1}^{R_0-1} \tau_{ij}.
\end{equation}
The truncation tolerances $\tau_{ij}$ for each of these sums must be chosen according to Eqn. \ref{eqn:trunc-error-bound} such that
\begin{equation}
    \sum_{i=1}^r \left( \norm{e_i}_2^2 + 2 \sigma_i \norm{e_i}_2  \right )  \le \tau_\textrm{TT}.
\end{equation}
A simple choice is to enforce
\begin{equation}
    \label{eqn:col-V-const-trunc-tol}
    \norm{e_i}_2^2 + 2 \sigma_i \norm{e_i}_2 \le \frac{\tau_\textrm{TT}}{r} \ , \quad \textrm{or equivalently} \quad \norm{e_i}_2 \le -\sigma_i + \sqrt{\sigma_i^2 + \frac{\tau_\textrm{TT}^2}{r^2}} \ , 
\end{equation}
which can be achieved by setting
\begin{equation}
    \label{eqn:TT-sum-const-trunc-tol}
    \tau_{ij} = \frac{1}{R_0 - 1}\left(-\sigma_i + \sqrt{\sigma_i^2 + \frac{\tau_\textrm{TT}^2}{r^2}}\right) \ . 
\end{equation}

Alternatively, one can select the tolerances $\tau_{ij}$ adaptively as terms are added to the sum. The errors $\norm{e_{ij}}_2$ can be computed exactly with minimal extra cost when performing each TT-SVD, and these values can be used to inform the subsequent truncation tolerances.
Start with $\tau_{i,0}$ as in Eqn. \ref{eqn:TT-sum-const-trunc-tol}, and set subsequent truncation tolerances $\tau_{ij}$ for $j = 1, 2, ..., R_0-1$ as
\begin{equation}
    \tau_{ij} = j \cdot \tau_{i,0} - \sum_{k < j} \norm{e_{ij}}_2 \ . 
\end{equation}
The idea of adaptive truncation tolerances can be applied at the level of the density matrix as well. Instead of the fixed error $\tau_\textrm{TT}/r$ per column of $V$ as in Eqn. \ref{eqn:col-V-const-trunc-tol}, we can enforce
\begin{equation}
    \norm{e_i}^2 + 2 \sigma_i \norm{e_i}_2 \ \le \ \tau_\textrm{TT}^{(i)} \ := \ i \cdot \frac{\tau_\textrm{TT}}{r} - \sum_{\ell < i} \left( \norm{e_\ell}_2^2 + 2 \sigma_\ell \norm{e_\ell}_2  \right ) ,
\end{equation}
and pick the truncation tolerances $\tau_{ij}$ as
\begin{equation}
    \tau_{ij} = \frac{j}{R_0-1} \left(-\sigma_i + \sqrt{\sigma_i^2 + (\tau_\textrm{TT}^{(i)})^2}\right) - \sum_{k < j} \norm{e_{ij}}_2  . 
\end{equation}
Algorithm \ref{SM-alg:iterative-TT-SVD} provides the pseudo-code for such this adaptive truncation scheme using TT-SVD.

\begin{algorithm}[t!]
\caption{TT-SVD for performing the linear combinations within \texttt{TT-Compress}}
\label{SM-alg:iterative-TT-SVD}
\begin{algorithmic}[1]
\Require Factor matrix $X \in \C^{N \times R_0}$ with columns $x_i \in \C^N$ represented in TT/MPS format; Coefficient matrix $V \in \C^{R_0 \times r}$, Singular values $\sigma_1,...,\sigma_r$, \ Error tolerance $\tau$
\vspace{0.35\baselineskip}
\Ensure Factor matrix $\tilde{X} \in \C^{N \times r}$ with $\lVert XX^\dagger - (XV)(XV)^\dagger\rVert_F \le \tau$
\vspace{0.35\baselineskip}
\Procedure{Linear-Combinations-via-TT-SVD}{$X, V, \sigma, \tau$}
    \vspace{0.5\baselineskip}
    \State $\tau_{\textrm{per-term}} = \tau \ / \ r$ 
    \Comment Expected truncation per column $\tilde{x}_i$ of $\tilde{X}$
    \State $\epsilon_\textrm{upper-bound} = 0$ 
    \Comment{Accumulator of error due to TT/MPS truncation}
    \For {$i = 1, 2, ..., r$}
        \State $\tau_\textrm{TT}^{(i)} = -\sigma_i + \left(\sigma_i^2 + (i \cdot \tau_\textrm{per-term} - \epsilon_\textrm{upper-bound})\right)^{1/2}$
        \vspace{0.25\baselineskip}
        \State $[\tilde{x}_i, \epsilon_i] = \texttt{TT-Sum}\left(\sum_{j=1}^R V_{ji} x_j, \ \tau_\textrm{TT}^{(i)}\right)$ 
        \Comment{$\epsilon_i = \norm{\tilde{x}_i - \sum_{j=1}^R V_{ji} x_j}_F \le \tau_\textrm{TT}^{(i)}$}
        \vspace{0.25\baselineskip}
        \State $\epsilon_\textrm{upper-bound} \gets \epsilon_\textrm{upper-bound} + 2 \sigma_i \epsilon_i + \epsilon_i^2 $
	\EndFor
\EndProcedure
\end{algorithmic}
\end{algorithm}

%
%
\section{Operator Splitting with Nearest-Neighbors Interactions}
\label{SM-sec:TEBD}

A standard approach to constructing the flow operator $\exp\{-i h H_\textrm{eff}\}$ is to use the operator splitting technique, commonly called \textit{time-evolving block decimation} (TEBD) in the tensor network community \cite{Vidal2004-TEBD, Verstraete2004-TEBD}.
This approach is most suitable for systems with nearest neighbor interactions between dimensions, e.g. coupling only between sites $i$ and $i+1$ of the MPS.
In general, one must partition the Hamiltonian as
\begin{equation}
    H_\textrm{eff} = \sum_{p=1}^{P} \sum_{m=1}^{M_p} H^{(p,m)} \quad \textrm{such that} \quad [H^{(p,m)}, H^{(p,m')}] = 0 .
\end{equation}
This is to say that each subset $\{H^{(p,m)} \ | \ m = 1, 2, ..., M_p\}$ contains terms that all commute with each other.
In such case, the flow operator with respect to each subset can be computed exactly as a product
\begin{equation}
    U_{h}^{(p)} := \exp\left\{-ih  \sum_{m=1}^{M_p} H^{(p,m)}\right\} = \prod_{m=1}^{M_p} \exp \left \{-i h H^{(m,p)} \right \}.
\end{equation}
The overall flow operator can then be approximated to two second order as
\begin{equation}
    U_{h} = \exp\{-ih H_\textrm{eff}\} = \left(\prod_{p=1}^P U_{h/2}^{(p)}\right)\left(\prod_{q=0}^{P-1} U_{h/2}^{(P-q)}\right) + O(h^3).
\end{equation}
One can similarly construct higher-order splitting schemes using products of the $U_{c h}^{(p)}$ but for different step sizes $c h$.
Of note, negative step sizes (e.g. $c < 0$) are necessary to build splitting schemes of order 4 or higher, so such schemes may unstable for dissipative systems like open quantum systems.
We did not observe this instability in our numerical experiments, however.

A crux of splitting schemes is ones ability to efficiently compute and multiply the sub-flow operators $\exp \left \{-i h H^{(m,p)} \right \}$.
This is often the case when working with quantum systems because terms in the Hamiltonian typically act only on one or two dimensions.

Consider a Hamiltonian with nearest neighbors interactions $H = \sum_{k=1}^{d-1} c_{k,k+1}$ where 
\begin{equation}
    c_{k,k+1} = I_{n_d} \otimes \cdots \otimes I_{n_{k+2}} \otimes C_{k,k+1} \otimes I_{n_{k-1}} \otimes \cdots \otimes I_{n_1} \ , 
\end{equation}
with $C_{k,k+1} \in \C^{n_k n_{k+1} \times n_k n_{n+1}}$. 
The action of the operator $c_{k,k+1}$ depends only on the $k$-th and $(k+1)$-st dimensions. The flow operator $\exp\{-i h \ c_{k,k+1}\}$ can therefore be written in terms of the smaller matrix $\exp\{-ih \ C_{k,k+1}\}$ as
\begin{equation}
    \exp\{ -ih \ c_{k,k+1} \} = I_{n_d} \otimes \cdots \otimes I_{n_{k+2}} \otimes \exp\{-ih \ C_{k,k+1}\} \otimes I_{n_{k-1}} \otimes \cdots \otimes I_{n_1} .
\end{equation}
For quantum problems, the matrices $C_{k,k+1}$ are typically small enough such that $\exp\{-ih \ C_{k,k+1}\}$ can be computed quickly and to machine precision. 
Moreover, the terms $c_k$ in $H$ can be partitioned into two sets of pairwise commuting operators by partitioning the sum over $k$ into odd and even terms, namely
\begin{equation}
    H = H_\textrm{odd} + H_\textrm{even} \ , \quad H_\textrm{odd} = \sum_{k \ \textrm{odd}} c_{k,k+1} \ , \quad H_\textrm{even} = \sum_{k \ \textrm{even}} c_{k,k+1}.
\end{equation}
When the number of dimensions $d$ is even, flow operators with respect to either the odd or even Hamiltonians can be written explicitly as
\begin{align}
    \label{eqn:tebd-odd}
    &\exp\{-i h H_\textrm{odd}\}
    \\
    &\phantom{ABC}= \exp\{-ih \ C_{n-1,n}\} \otimes \cdots \otimes \exp\{-ih \ C_{3,4}\} \otimes \exp\{-ih \ C_{1,2}\},
    \nonumber \\[0.5em]
    \label{eqn:tebd-even}
    &\exp\{-i h H_\textrm{even}\} 
    \\
    &\phantom{ABC} = I_{n_d} \otimes \exp\{-ih \ C_{n-2,n-1}\} \otimes \cdots \otimes \exp\{-ih \ C_{2,3}\} \ \otimes I_{n_1} \nonumber .
\end{align}
and similar expressions exist for odd $d$. In either case, both $\exp\{-i h H_\textrm{odd}\}$ and $\exp\{-i h H_\textrm{even}\}$ can easily be written MPOs in terms of each $\exp\{-i h \ C_{k,k+1}\}$, represented as an order-2 MPO
\begin{equation}
    \big[\exp\{-i h \ C_{k,k+1}\}\big](i_k,i_{k+1};j_k,j_{k+1}) \ = \ \sum_{\sigma_k=1}^{b_k} L_k(i_k,j_k,\sigma_k) R_k(\sigma_k, i_{k+1},j_{k+1}).
\end{equation}
For instance, the odd flow operator in Eqn. \ref{eqn:tebd-odd} will be the MPO with cores $L_1, R_1, L_3$, $R_3, ..., L_{n-1}, R_{n-1}$. It will have bond dimensions $\{r_1, 1, r_3, 1, ..., 1, r_{d-1}\}$. On the other hand, the even flow operator in Eqn. \ref{eqn:tebd-even} will have cores $I_{n_1}$, $L_2$, $R_2$, $L_4$, $R_4, ..., L_{n-2}$, $R_{n-2}$, $I_{n}$, which has bond dimensions $\{1,r_2,1,r_4,...,1,r_{n-2},1\}$.

Once the odd and even flow operators have been constructed as MPOs, they are used within any existing splitting scheme.
For schemes of order 2 or more, MPO truncation after each MPO-MPO multiplication may be needed to reduce MPO bond dimensions.
The truncation tolerance should be $\tau = O(h^{p+1})$, where $p$ is the order of the method, to observe proper convergence.

%
%
\section{Mock-Circuit Simulation --- Randomized rounding of MPO - MPS products}

For our second experiment (c.f. Sec. \ref{sec:mock-quantum-circuit-sim}), the flow operators can have MPO bond dimension upwards of 50, so compression of the MPO-MPS products $y = U_h x$ becomes expensive even when the MPS $x$ has small bond dimension. 
Indeed, if $x$ has only bond dimension $\chi = 10$, the cores of the uncompressed $U_h x$ already have size $\approx 500 \times 2 \times 500$, and truncation of one product alone can take multiple seconds.
A better approach to this compression is to use randomized MPS rounding techniques introduced by \cite{Daas2021-Randomized-TT-Rounding}. We find that this choice drastically reduces the cost of compressing $U_h x$, often being upwards of 50x faster with errors nearly identical to the deterministic rounding procedure.

%
%
\section{Mock-Circuit Simulation --- Additional Data}

\begin{figure}[t!]
    \centering
    \includegraphics[width=\textwidth]{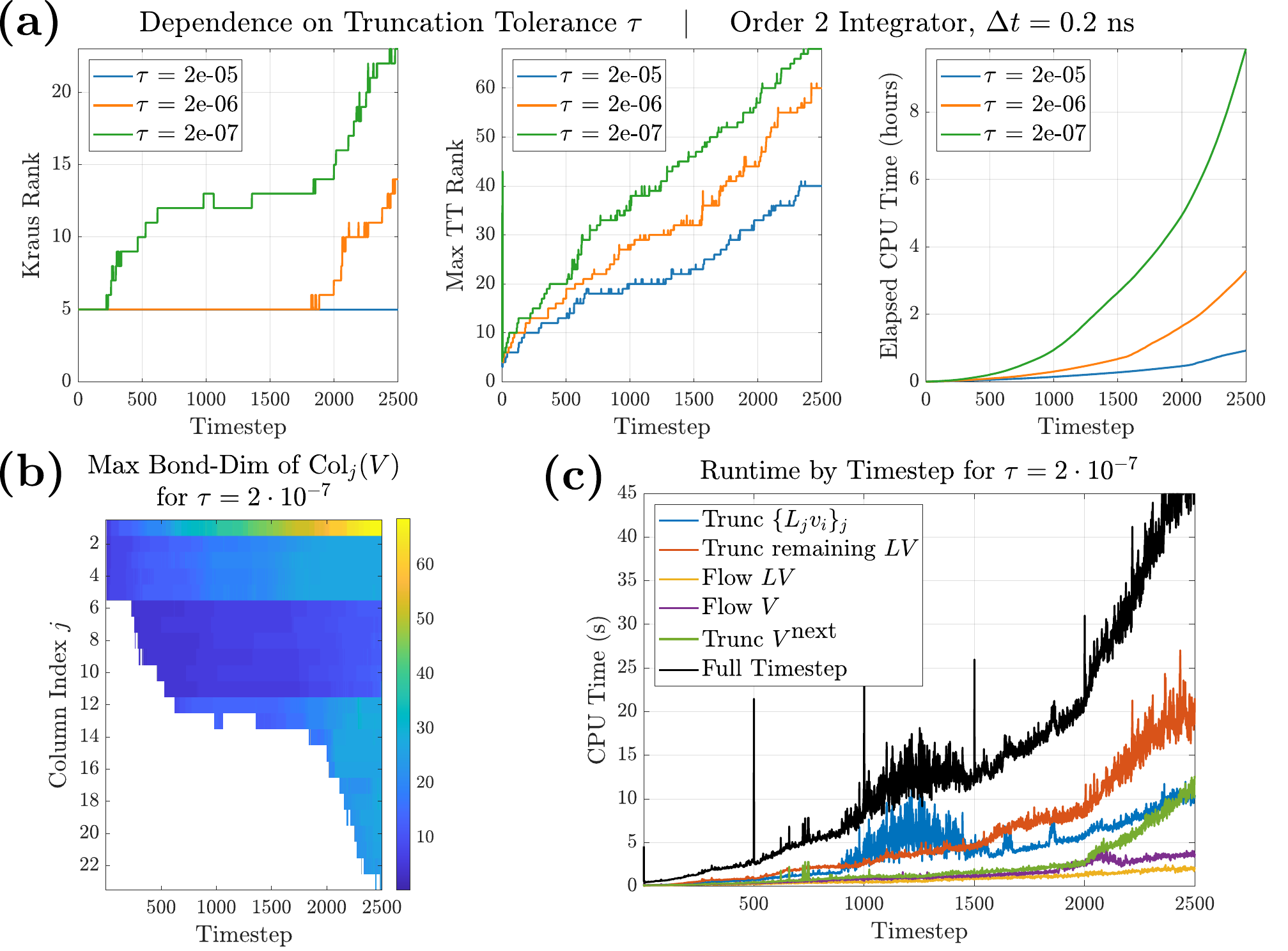}
    \caption{Run statistics of the heavy-hex lattice mock quantum circuit simulation. \textbf{(a)} shows the density matrix rank, maximum MPS bond dimension, and elapsed runtime by timestep when using three different rounding tolerances.
    \textbf{(b)} plots the maximum bond dimension by column of the factor matrix V for the $\tau = 2 \cdot 10^{-7}$ simulation, and \textbf{(c)} shows the runtime breakdown per timestep for different subroutines of the CPTP scheme.
    }
    \label{SM-fig:heavy-hex-run-stats}
\end{figure}

Fig. \ref{SM-fig:heavy-hex-run-stats} shows the dependence of the density matrix rank, maximum MPS bond dimension, and simulation runtime on the rounding tolerance $\tau$. 
All three metric increase as the tolerance decreases, with the smallest tested tolerance $\tau = 2 \cdot 10^{-7}$ using far more computational resources than the rest.

Fig. \ref{SM-fig:heavy-hex-run-stats}\textbf{(b)} shows the maximum bond dimensions of the columns $v_j(t)$ of the factor matrix $V(t)$. As before, each row in the heatmap delineates a column $v_j(t)$ of $V(t)$, with the color white at cell $(j,n)$ in the heatmap means that $V(t_n)$ had $r(t) < j$ columns.
Most columns have maximum bond at most 30, though the primary component has bond dimension that exceeds 60 by the end of the simulation.
The use of randomized rounding to compress MPO-MPS products is particularly important for these vectors, as otherwise the flow with respect to the effective Hamiltonian would be the dominating cost of the algorithm. In Fig. \ref{SM-fig:heavy-hex-run-stats}\textbf{(c)}, we see that performing said flow, while expensive, is small compared to the low-rank truncations of the density matrix.
Indeed, by the end of the simulation, the integrator spends 30+ seconds per timestep on these truncations, even when using randomized methods when computing linear combinations.

%
%
\section{Qudit-Resonator Chains --- Operator Splitting}
\label{supp:qudit-resonator-tebd}

In this section, we describe the operator splitting approach we took to building approximate MPO representations for the flow operators of the qubit-resonator chain simulations. Recall, for this problem, the Hamiltonian decomposes as $H(t) = H_d + H_c(t)$. The first term, called the \textit{drift Hamiltonian}, takes the form
\begin{equation}
    H_d = - \frac{1}{2}\sum_{q} \xi_q a_q^\dagger a_q^\dagger a_q a_q - \sum_{\langle p, q \rangle} \xi_{pq} a_p^\dagger a_p a_q^\dagger a_q \ , 
\end{equation}
where $a_q$ denotes the lowering operator for subsystem $q$
\begin{equation}
    a_q = I_{n_d} \otimes \cdots \otimes I_{n_{q+1}} \otimes A_q \otimes I_{n_{q-1}} \otimes \cdots \otimes I_{n_1} \ , 
\end{equation}
with
\begin{equation}
    A_q = \sum_{j=1}^{n_q-1} \sqrt{j} \ \ketbra{j}{j+1} \ \in \ \mathbb{R}^{n_q \times n_q}.
\end{equation}
$\xi_q$ is the self-Kerr coefficient of subsystem $q$, and $\xi_{pq}$ is the cross-Kerr coupling strength between subsystems $p$ and $q$.
The sum introducing the coupling is taken over pairs $\langle p, q \rangle$ determined by an underlying device layout.

The control Hamiltonian $H_c$ introduces off-diagonal elements. It is defined as
\begin{equation}
    H_c(t) = \sum_{q=1}^d \big( d_q(t) a_q + \bar{d}_q(t) a_q^\dagger\big).
\end{equation}
Each $d_q(t)$ is a complex-valued function denoting the control signal to transmon $q$. For our experiments, these functions are piecewise constant, with value changing every $T_\textrm{signal}$ nanoseconds.

The Hamiltonian remains constant with each interval $\big[(k-1) T_\textrm{signal}, kT_\textrm{signal} \big)$ for $k = 1, ..., T_\textrm{gate}/T_\textrm{signal}$. Letting $H^{(k)}$ denote the Hamiltonian during the $k$-th time interval, we build an approximation of the flow operator $U^{(k)} := \exp\{-i h H^{(k)}\}$ via operator splitting.
The qudit-qudit cross-Kerr terms introduce interactions between MPS sites a distance 2 apart, so a bit more care is needed for this construction as compared to the spin-chain example.
Here, we partition $H^{(k)}$ into a three parts as
\begin{align}
    H^{(k)} &= \overbrace{H_1 + H_2}^{H_d} + \overbrace{H_3^{(k)}}^{H_c(t)} ,
    \\[0.5em]
    H_1 &= \sum_{r \equiv 0 \bmod{2}} H_{CK}^{(r)} ,
    \\[0.5em]
    H_2 &= \sum_{r \equiv 1 \bmod{2}} (H_{SK}^{(r)} + H_{CK}^{(r)}),
    \\[0.5em]
    H_3^{(k)} &= H_c(t_k) \ , \enskip t_k =(k-1)T_\textrm{signal}.
\end{align}
Here, $H_{CK}^{(r)}$ denotes the cross-Kerr coupling between sites $2r-1, 2r$, and $2r+1$, namely
\begin{align}
    H_{CK}^{(r)} = \
    -\xi_{q-1,q} (a_{q-1}^\dagger a_{q-1}) (a_{q}^\dagger a_{q}) 
    & - \ \xi_{q,q+1} (a_{q}^\dagger a_{q}) (a_{q+1}^\dagger a_{q+1})   & (\textrm{where} \ q = 2r) \nonumber
    \\ 
    & - \ \xi_{q-1,q+1} (a_{q-1}^\dagger a_{q-1}) (a_{q+1}^\dagger a_{q+1}).  \nonumber
\end{align} 
The self-kerr terms $H_{SK}^{(r)}$ is defined similarly as
\begin{align}
    H_{SK}^{(r)} = -\sum_{q=2r-1}^{2r+1} \frac{\xi_q}{2} a_q^\dagger a_q^\dagger a_q a_q.
\end{align} 
The matrix exponential of each $H_{CK}^{(r)}$ is cheap to compute because these operators only act on three adjacent sites of the MPS. In particular, building $\exp\{-i h H_{CK}^{(r)}\}$ amounts to computing the flow operator with respect to the 3-site Hamiltonian,
\begin{align*}
    H_{CK}^{(r,3\textrm{-site})} =
    & -\xi_{2r-1,2r} \ I_{n_\textrm{qudit}} \otimes (A_{n_\textrm{res}}^\dagger A_{n_\textrm{res}}) \otimes (A_{n_\textrm{qudit}}^\dagger A_{n_\textrm{qudit}}) 
    \\
    &- \ \xi_{2r,2r+1} \ (A_{n_\textrm{qudit}}^\dagger A_{n_\textrm{qudit}}) \otimes (A_{n_\textrm{res}}^\dagger A_{n_\textrm{res}}) \otimes I_{n_\textrm{qudit}}
    \\ 
    &- \ \xi_{2r-1,2r+1} \ (A_{n_\textrm{qudit}}^\dagger A_{n_\textrm{qudit}}) \otimes I_{n_\textrm{res}} \otimes (A_{n_\textrm{qudit}}^\dagger A_{n_\textrm{qudit}}),
\end{align*}
of size $n_\textrm{qudit}^2 n_\textrm{res} \times n_\textrm{qudit}^2 n_\textrm{res}$ and then decomposing the resulting operator into an MPO with mode sizes $(n_\textrm{qudit},n_\textrm{res},n_\textrm{qudit})$. 
The cores of this 3-site smaller MPO are used as the MPO cores for sites $2r-1,2r,2r+1$ of the MPO representation of $\exp\{-i h H_{CK}^{(r)}\}$, with all other cores being identity operators of the appropriate sizes.
The matrix exponential of $H_{SK}^{(r)} + H_{CK}^{(r)}$ can be computed in a similar manner.

Any two cross-Kerr Hamiltonians $H_{CK}^{(r)}$ and $H_{CK}^{(r')}$ commute with each other so long as $|r - r'| > 1 $. As such,
\begin{align}
\exp\{-i h H_1 \} &= \prod_{r \equiv 0 \bmod{2}} \exp\{-i h H_{CK}^{(r)} \},
\\[0.5em]
\exp\{-i h H_2 \} &= \prod_{r \equiv 1 \bmod{2}} \exp\bigg\{-i h \big(H_{SK}^{(r)} +H_{CK}^{(r)}\bigg) \bigg \}.
\end{align}
The terms in these products need not be explicitly formed as full MPOs. For instance, when building $\exp\{-i h H_1\}$ as an MPO, one only needs to compute the non-identity cores of each $\exp\{-i h H_{CK}^{(r)}\}$, e.g. those at sites $2r-1, 2r,$ and $2r+1$, as these are the corresponding cores in the MPO representation of $\exp\{-i h H_1 \}$. 

The final term $H_3^{(k)}$ in our splitting is already a sum of pairwise computing terms --- the control signals which act independent on each site. It's matrix exponential is simply the rank-1 MPO
\begin{align}
    \exp\{-i h H_3^{(k)} \} &= \prod_{q=1}^d \exp\left\{-i h \left(d_q(t_k) a_q + \bar{d}_q(t_k)a_q^\dagger\right)\right\} \nonumber
    \\[0.5em]
    &= \bigotimes_{q=d,...,2,1} \exp\left\{-i h \left(d_q(t_k) A_{n_q} + \bar{d}_q(t_k)A_{n_q}^\dagger\right)\right\} \quad (A_{n_q} \in \mathbb{C}^{n_q \times n_q}).
\end{align}

Now that we can build MPO representations for the matrix exponentials of each term in our splitting, we can obtain a second-order approximation for the full flow operator via Strang splitting
\begin{align}
    U^{(k)} &:= \exp\{-i h H^{(k)}\} = \exp\left\{-i \frac{h}{2} H_3^{(k)} \right\} \ U_{d} \ \exp\left\{-i \frac{h}{2} H_3^{(k)} \right\} + O(h^3),
    \\[0.5em]
    U_{d} &:= \exp\{-i h H_d\} = \exp\left\{-i \frac{h}{2} H_1 \right\} \ \exp\left\{-i h H_2 \right\} \ \exp\left\{-i \frac{h}{2} H_1 \right\} + O(h^3).
\end{align}
MPO compression of the device flow operator $U_d$ is generally necessary to maintain low MPO ranks. This operator only needs to be constructed once, however, as the device parameters do not change in time.
The only time-variance comes from the control term $H_3^{(k)}$, whose matrix exponential multiplies $U_{d}$ to the left and right appears in the equation for $U^{(k)}$. Luckily $\exp\{-ih H_3^{(k)}\}$ is a rank-1 MPO, so computing $U^{(k)}$ from a pre-computed $U_d$ requires no additional MPO compression.

\end{document}